\documentclass[11pt]{article}
\usepackage{amsmath}
\usepackage{amssymb}
\usepackage{latexsym}
\usepackage{dsfont}
\def\disp{\displaystyle}

\def\Limsup{\mathop{{\rm Lim}\,{\rm sup}}}

\def\tto{\;{\lower 1pt \hbox{$\rightarrow$}}\kern -10pt
\hbox{\raise 2pt \hbox{$\rightarrow$}}\;}

\def\Hat{\widehat}

\def\Bar{\overline}
\def\ra{\rangle}
\def\la{\langle}
\def\ve{\varepsilon}
\def\B{I\!\!B}
\def\h{\hfill\Box}
\def\R{\mathbb{R}}
\def\N{\mathbb{N}}

\def\ox{\bar{x}}
\def\oy{\bar{y}}
\def\oz{\bar{z}}

\def\co{\mbox{\rm co}\,}
\def\gph{\mbox{\rm gph}\,}
\def\epi{\mbox{\rm epi}\,}

\def\dom{\mbox{\rm dom}\,}

\def\cl{\mbox{\rm cl}\,}
\def\dist{\mbox{\rm dist}\,}
\def\cone{\mbox{\rm cone}\,}

\def\bd{\mbox{\rm bd}\,}
\def\inter{\mbox{\rm int}\,}

\def\substack#1#2{{\scriptstyle{#1}\atop\scriptstyle{#2}}}

\def\h{\hfill\square}
\def\dn{\downarrow}
\def\O{\Omega}

\def\ph{\varphi}
\def\emp{\emptyset}
\def\st{\stackrel}
\def\oR{\Bar{\R}}
\def\lm{\lambda}
\def\Lm{\Lambda}

\def\Th{\Theta}
\def\GG{\Gamma}
\def \N{I\!\!N}

\begin{document}
\begin{center}
{\bf TANGENTIAL EXTREMAL PRINCIPLES FOR FINITE AND INFINITE
SYSTEMS OF SETS, II: APPLICATIONS TO SEMI-INFINITE AND
MULTIOBJECTIVE OPTIMIZATION}\footnote{This research was partially
supported by the US National Science Foundation under grants
DMS-0603846 and DMS-1007132 and by the Australian Research Council
under grant DP-12092508.}\\[2ex]
BORIS S. MORDUKHOVICH\footnote{Department of Mathematics, Wayne
State University, Detroit, MI 48202, USA. Email:
boris@math.wayne.edu.} and HUNG M. PHAN\footnote{Department of
Mathematics, Wayne State University, Detroit, MI 48202, USA.
Email: pmhung@wayne.edu.}
\end{center}
\newtheorem{Theorem}{Theorem}[section]
\newtheorem{Proposition}[Theorem]{Proposition}
\newtheorem{Lemma}[Theorem]{Lemma}
\newtheorem{Corollary}[Theorem]{Corollary}
\newtheorem{Definition}[Theorem]{Definition}
\newtheorem{Remark}[Theorem]{Remark}
\newtheorem{Example}[Theorem]{Example}
\renewcommand{\theequation}{\thesection.\arabic{equation}}
\normalsize
\setcounter{equation}{0}
\setcounter{section}{0}
\newcounter{count}
\newenvironment{lista}{\begin{list}{\rm(\alph{count})}{\usecounter{count}
\setlength{\itemindent}{-0.3cm}}}{\end{list}}
\newenvironment{listi}{\begin{list}{\rm(\textit{\roman{count}})}{\usecounter{count}
\setlength{\itemindent}{-0.3cm}}}{\end{list}}
\newenvironment{listH}{\begin{list}{\bf{(H\arabic{count})}}{\usecounter{count}
\setlength{\itemindent}{0cm}}}{\end{list}} \small{{\bf Abstract.}
This paper contains selected applications of the new tangential
extremal principles and related results developed in
\cite{MorPh10a} to calculus rules for infinite intersections of
sets and optimality conditions for problems of semi-infinite
programming and multiobjective optimization with countable
constraints.\vspace*{0.1cm}

\vspace*{0.1cm}

}
\normalsize

\section{Introduction}
\setcounter{equation}{0}

Variational analysis is based on variational principles and
techniques, which are largely inspired and motivated by
applications to constrained optimization and related problems.
Extremal principles for systems of sets can be treated as
variational principles in a geometric framework while playing a
crucial role in the core variational theory and numerous
applications; see, in particular, the books
\cite{Borwein-Zhu-TVA,m-book1,m-book2,Rockafellar-Wets-VA,s} and
the references therein.

In \cite{MorPh10a}, we developed new tangential extremal
principles that concerned, for the first time in the literature,
{\em countable} systems of sets. Our main motivation came from
possible applications to problems of semi-infinite optimization
with a countable number of constraints. It has been well
recognized in optimization theory and its applications that
problems of this type are significantly more difficult in
comparison with conventional problems of semi-infinite
optimization dealing with parameterized constraints over compact
index sets; see, e.g., \cite{gl}.

This paper mainly addresses selected applications of the
tangential extremal principles and their consequences in
\cite{MorPh10a} to various problems of semi-infinite optimization
with countable constraints, particularly including those which
naturally arise in semi-infinite programming and multiobjective
optimization. To deal with such problems, we develop new calculus
rules for tangent and normal cones to countable intersections of
sets. These calculus results are certainly of their own interest
being strongly used in the subsequent applications. To simplify
the presentation, we confine ourselves to problems formulated in
finite-dimensional spaces. At the same time, the initial data
involved may be nonsmooth and nonconvex, and we strongly employ
appropriate constructions of generalized differentiation in
variational analysis.\vspace*{0.05in}

The rest of the paper is organized as follows. Section~2 contains
some preliminaries of variational analysis and also recall two
major results from \cite{MorPh10a} largely used in the sequel.

Section~3 is devoted to calculus rules for tangent and normal
cones to countable intersections of nonconvex sets and the
corresponding qualification conditions. A special attention is
paid in this section to a countable nonconvex version of the
so-called ``conical hull intersection property" (CHIP) developed
earlier for finite intersections of convex sets and successfully
used in convex optimization, approximation theory, etc. We
establish verifiable sufficient conditions for the nonconvex CHIP
and employ this property and other qualification conditions to
derive new calculus rules for generalized normals to infinite
intersections of nonconvex sets in finite dimensions.

Section~4 presents a number of applications of the results from
\cite{MorPh10a} and the from the preceding section to deriving
necessary optimality conditions in various problems of
semi-infinite programming with geometric, operator, and functional
constraints. We obtain optimality conditions of different types
under appropriate constraint qualifications and compare the
optimality and qualification conditions obtained with those known
before in convex and nonconvex settings.

Finally, Section~5 concerns applications of our major tangential
extremal principle and the related calculus rules to various
problems of multiobjective optimization including those with
set-valued objectives. Besides paying the main attention to
multiobjective problems with countable constraints, we introduce
and develop there some notions and results, which seem to be of
their own interest for the general theory of multiobjective
optimization and its subsequent applications.\vspace*{0.05in}

The notation and terminology of the paper are basically standard
in variational analysis and generalized differentiation; cf.\
\cite{MorPh10a} and the books on variational analysis mentioned
above. Recall that $\N:=\{1,2,\ldots\}$, that $\B$ denotes the
closed unit ball in $\R^n$, and that
\begin{equation}\label{eq:LimsupPK}
\begin{aligned}
\Limsup_{x\to\ox}F(x):=\Big\{y\in\R^m\,\Big|\,&\exists\mbox{ sequences }\;x_k\to\ox\;\mbox{ and }\;y_k\to y\\
&\mbox{with }\;y_k\in F(x_k)\;\mbox{ for all }\;k\in\N\Big\}
\end{aligned}
\end{equation}
stands for the (sequential) {\em Painlev\'{e}-Kuratowski
upper/outer limit} of a set-valued mapping $F\colon\R^n\tto\R^m$
at a point $\ox\in\dom F:=\{x\in\R^n|\;F(x)\ne\emp\}$ of its
domain.

\section{Preliminaries from Variational Analysis}
\setcounter{equation}{0}

Our main references for the brief overview of this section are
\cite{m-book1,MorPh10a,Rockafellar-Wets-VA}, where the reader can
find proofs, discussions, and commentaries.

Given a set $\O\subset\R^n$ locally closed around a point
$\ox\in\O$, we use in this paper the (only one) notion of the {\em
tangent cone} $T(\ox;\O)$ given by
\begin{equation}\label{tc}
T(\ox;\O):=\Limsup_{t\dn0}\frac{\O-\ox}{t},
\end{equation}
 which is also known as the Bouligand-Severi contingent cone to $\O$ at $\ox$.
 The {\em normal cone} $N(\ox;\O)$ to $\O$ at $\ox$ is defined by
the outer limit \eqref{eq:LimsupPK} as
\begin{equation}\label{nc}
N(\ox;\O):=\Limsup_{x\to\ox}\big[\cone\big(x-\Pi(x;\O)\big)\big]
\end{equation}
via the Euclidean projection
$\Pi(x;\O):=\{w\in\O|\;\|x-w\|=\mbox{dist}(x;\O)\}$ to $\O$ at
$x\in\O$ and is known under that names of the
Mordukhovich/basic/limiting normal cone to closed subsets of
finite-dimensional spaces. Our basic normal cone \eqref{nc} is
often nonconvex while admitting the following outer limiting
representation:
\begin{equation*}
N(\ox;\O)=\Limsup_{x\st{\O}{\to}\ox}\Hat N(x;\O)
\end{equation*}
via the convex collections of {\em Fr\'echet normals} to $\O$ at $x\in\O$ given by
\begin{equation}\label{eq:e-nor}
\Hat N(x;\O):=\left\{x^*\in\R^n\Big|\;\limsup_{u\st{\O}{\to}x}\frac{\la
x^*,u-x\ra}{\|u-x\|}\le 0\right\},
\end{equation}
where $u\st{\O}{\to}x$ means that $u\to x$ with $u\in\O$. Note
that $\Hat N(\ox;\O)$, known also as the prenormal or regular
normal cone, is dually generated by the (generally nonconvex)
tangent cone \eqref{tc} as
\begin{equation}\label{du}
\Hat N(\ox;\O)=T^*(\ox;\O):=\big\{x^*\in\R^n\big|\;\la x^*,v\ra\le 0\;\mbox{ for all }\;v\in T(\ox;\O)\big\}.
\end{equation}
For convex sets $\O$ all the constructions
\eqref{tc}--\eqref{eq:e-nor} reduce to the corresponding tangent
and normal cones of convex analysis, while only the basic normal
cone \eqref{nc} enjoys comprehensive calculus rules (full
calculus) in general nonconvex settings; see
\cite{m-book1,Rockafellar-Wets-VA} and their references. Note the
following remarkable fact relating the tangent and normal cones to
arbitrary closed sets $\O\subset\R^n$ (see
\cite[Theorem~6.27]{Rockafellar-Wets-VA} and
\cite[Corollary~6.5]{MorPh10a}):
\begin{equation}
\label{eq:N(Lm)-N(O)}
N\big(0;T(\ox;\O)\big)\subset N(\ox;\O).
\end{equation}

Given further a set-valued mapping $F\colon\R^n\tto\R^m$ with the graph
\begin{equation*}
\gph F:=\big\{(x,y)\in\R^n\times\R^m\big|\;y\in F(x)\big\},
\end{equation*}
we define the {\em coderivative} of $F$ at $(\ox,\oy)\in\gph F$ via the normal cone \eqref{nc} by
\begin{equation}\label{eq:Coder}
D^*F(\ox,\oy)(y^*):=\big\{x^*\in\R^n\big|\;(x^*,-y^*)\in
N\big((\ox,\oy);\gph F\big)\big\},\quad y^*\in\R^m,
\end{equation}
where $\oy=f(\ox)$ is omitted if $F=f\colon\R^n\to\R^m$ is
single-valued. Observe that the coderivative \eqref{eq:Coder} is a
positively homogeneous mapping $D^*F(\ox,\oy)\colon\R^m\tto\R^n$,
which reduces to the single-valued adjoint derivative operator
\begin{equation}\label{smooth}
D^*f(\ox)(y^*)=\big\{\nabla f(\ox)^*y^*\big\}\;\mbox{ for all }\;y^*\in\R^m
\end{equation}
if $f$ is {\em strictly differentiable} at $\ox$ in the sense that
\begin{equation*}
\lim_\substack{x\to\ox}{u\to\ox}\frac{f(x)-f(u)-\la\nabla f(\ox),x-u\ra}{\|x-u\|}=0;
\end{equation*}
the latter is automatic if $f$ when $C^1$ around $\ox$.

Given finally an extended-real-valued function
$\ph\colon\R^n\to\oR:=[-\infty,\infty]$ finite at $\ox$, we define
its {\em basic subdifferential} at $\ox$ by
\begin{equation}\label{sub}
\partial\ph(\ox):=\big\{x^*\in\R^n\big|\;(x^*,-1)\in N\big((\ox,\ph(\ox));\epi\ph\big)\big\}
\end{equation}
via the normal cone \eqref{nc} to the epigraph
$\epi\ph:=\{(x,\mu)\in\R^{n+1}|\;\mu\ge\ph(x)\}$. The
subdifferential \eqref{sub} can be equivalently represented as the
outer limit
\begin{equation*}
\partial\ph(\ox)=\Limsup_{x\st{\ph}{\to}\ox}\Hat\partial\ph(x),
\end{equation*}
with $x\st{\ph}{\to}\ox$ indicating that $x\to\ox$ and $\ph(x)\to\ph(\ox)$, of the Fr\'echet-like construction
\begin{equation}\label{fr}
\Hat\partial\ph(x):=\big\{x^*\in\R^n\Big|\;\liminf_{u\to x}\frac{\ph(u)-\ph(x)-\la x^*,u-x\ra}{\|u-x\|}\ge 0\Big\}.
\end{equation}

To conclude this section, we recall the concept of tangential
extremality for countable set systems introduced in
\cite{MorPh10a} and formulate two major results obtained therein,
which are largely used in what follows. A set system
$\{\O_i\}_{i\in\N}\subset\R^n$ is {\em tangentially extremal} at
$\ox\in\bigcap_{i=1}^\infty\O_i$ if there is a bounded sequence
$\{a_i\}_{i\in\N}\subset\R^n$ such that
\begin{equation}\label{ext}
\bigcap_{i=1}^{\infty}\big[T(\ox;\O_i)-a_i\big]=\emp.
\end{equation}
\begin{Theorem}\label{Thm:TEPinf} {\bf (tangential extremal principle).} Let a countable system
$\{\O_i\}_{i\in\N}$ of closed sets in $\R^n$ be
tangentially extremal at $\ox\in\bigcap_{i=1}^\infty$. Assume that
\begin{equation}\label{ext1}
\bigcap_{i=1}^\infty\big[T(x;\O_i)\big]=\{0\}.
\end{equation}
Then there are  normal vectors
\begin{equation}\label{ext2}
x^*_i\in N\big(0;T(\ox;\O_i)\big)\subset N(\ox;\O_i)\;\mbox{ for all }\;i=1,2,\ldots
\end{equation}
satisfying the following extremality conditions:
\begin{equation}\label{eq:EPinf}
\sum_{i=1}^{\infty}\frac{1}{2^i}x^*_i=0\quad\mbox{and}\quad\sum_{i=1}^{\infty}\frac{1}{2^i}\|x^*_i\|^2=1.
\end{equation}
\end{Theorem}

The next result from \cite{MorPh10a} is based on the tangential extremal principle.

\begin{Theorem}\label{Thm:FnorLm} {\bf (representation of Fr\'echet normals to countable cone intersections).}
Let $\{\Lm_i\}_{i\in\N}$ be a countable system of closed cones in
$\R^n$ satisfying the conic qualification condition
\begin{equation}\label{eq:QC-N(Lm)}
\left[\sum_{i=1}^{\infty}x^*_i=0,\;x^*_i\in N(0;\Lm_i)\right]\Longrightarrow\Big[x^*_i=0\;\mbox{ for all}\;i\in\N\Big].
\end{equation}
Denoting the cone intersection by
$\Lm:=\bigcap_{i=1}^\infty\Lm_i$, we have the following
representation of Fr\'echet normals to $\Lm$ at the origin:
\begin{equation}\label{fr-rep}
\Hat N(0;\Lm)\subset\cl\Big\{\sum_{i\in I}x^*_i\Big|\;x^*_i\in N(0;\Lm_i),\;I\in{\cal L}\Big\},
\end{equation}
where ${\cal L}$ is the collection of all the finite subsets of $\N$.
\end{Theorem}

\section{Tangents and Normals to Infinite Intersections of Sets}
\setcounter{equation}{0}

The main purpose of this section is to derive calculus rules for
representing generalized normals to countable intersections of
arbitrary closed sets under appropriate qualification conditions.
Besides employing the tangential extremal principle, one of the
major ingredients in our approach is relating calculus rules for
generalized normals to countable set intersections with the
so-called ``conical hull intersection property" defined in terms
of tangents to sets, which was intensively studied and applied in
the literature for the case of finite intersections of convex
sets; see, e.g.,
\cite{Bauschke-Borwein-Li99,Chui-Deutsch-Ward90,Deutsch-Li-Ward97,Ernst-Thera07,Li-Ng-Pong07}
and the references therein. In what follows, we keep the
terminology of convex analysis (that goes back probably to
\cite{Chui-Deutsch-Ward90}) replacing the tangent and normal cones
therein by the nonconvex extension \eqref{tc} and \eqref{nc}.

\begin{Definition}\label{chip} {\bf (CHIP for countable intersections).} A set system $\{\O_i\}_{i\in\N}$
in $\R^n$ is said to have the {\sc conical hull intersection
property} $($CHIP$)$ at $\ox\in\bigcap_{i=1}^\infty\O_i$ if
\begin{equation}\label{eq:CHIP}
T\Big(\ox;\bigcap_{i=1}^\infty\O_i\Big)=\bigcap_{i=1}^\infty T(\ox;\O_i).
\end{equation}
\end{Definition}

In convex analysis and its applications the CHIP is often related
to the so-called ``strong CHIP" for finite set intersections
expressed via the normal cone to the convex sets in question.
Following this terminology in the case of infinite intersections
of nonconvex sets, we say that a countable system of sets
$\{\O_i\}_{i\in\N}$ has the {\em strong conical hull intersection
property} (or the {\em strong CHIP}) at
$\ox\in\bigcap_{i=1}^\infty\O_i$ if
\begin{equation}\label{eq:sCHIP}
N\Big(\ox;\bigcap_{i=1}^\infty\O_i\Big)=\Big\{\sum_{i\in I}x^*_i\Big|\;x^*_i\in N(\ox;\O_i),\;I\in{\cal L}\Big\}.
\end{equation}
When all the sets $\O_i$ as $i\in\N$ are {\em convex} in
\eqref{eq:sCHIP}, the strong CHIP of the system
$\{\O_i\}_{i\in\N}$ can be equivalently written in the form
\begin{equation}\label{eq:sCHIP1}
N\Big(\ox;\bigcap_{i=1}^\infty\O_i\Big)=\co\bigcup_{i=1}^\infty N(\ox;\O_i).
\end{equation}
We say that a countable set system $\{\O_i\}_{i\in\N}$ has the
{\em asymptotic strong CHIP} at $\ox\in\bigcap_{i=1}^\infty\O_i$
if the latter representation is replaced by
\begin{equation}\label{asym}
N\Big(\ox;\bigcap_{i=1}^\infty\O_i\Big)=\cl\co\bigcup_{i=1}^\infty N(\ox;\O_i).
\end{equation}

The next result shows the {\em equivalence} between the CHIP and
the asymptotic strong CHIP for intersections of convex sets in
finite dimensions. It follows from the proof that this equivalence
holds for {\em arbitrary} intersections of convex sets, not only
for countable ones studies in this paper.
\begin{Theorem}\label{Thm:CHIPchar} {\bf (characterization of CHIP for intersections of convex sets).}
Let $\{\O_i\}_{i\in\N}$ be a system of convex sets in $\R^n$, and let $\ox\in\bigcap_{i=1}^\infty\O_i$.
The following are equivalent:

{\bf (a)} The system $\{\O_i\}_{i\in\N}$ has the CHIP at $\ox$.

{\bf (b)} The system $\{\O_i\}_{i\in\N}$ has the asymptotic strong CHIP at $\ox$.\\[1ex]
In particular, the strong CHIP implies the CHIP but not vice versa.
\end{Theorem}
{\bf Proof}. Observe first that for convex sets in finite dimensions, in addition to the duality property
\eqref{du} with $\Hat N(\ox;\O)$ replaced by
$N(\ox;\O)$, we have the reverse duality representation
\begin{equation}\label{du1}
T(\ox;\O)=N^*(\ox;\O):=\big\{v\in\R^n\big|\;\la x^*,v\ra\le 0\;\mbox{ for all }\;x^*\in
N(\ox;\O)\big\}.
\end{equation}
Let us now justify the equality
\begin{equation}\label{eq:CHIPchar1}
\Big(\bigcap_{i=1}^\infty T(\ox;\O_i)\Big)^*=\cl\co\bigcup_{i=1}^\infty N(\ox;\O_i).
\end{equation}
The inclusion ``$\supset$" follows from \eqref{du} by the observation
\begin{equation*}
N(\ox;\O_i)=T^*(\ox;\O_i)\subset\Big(\bigcap_{i=1}^\infty T(\ox;\O_i)\Big)^*
\end{equation*}
due the closedness and convexity of the polar set on the right-hand side of the latter inclusion.

To prove the opposite inclusion ``$\subset$" in \eqref{eq:CHIPchar1}, pick some
$x^*\not\in\cl\co\bigcup_{i=1}^\infty N(\ox;\O_i)$. Then the classical separation theorem for
convex sets ensures the existence of a vector $v\in\R^n$ such that
\begin{equation}\label{v}
\la x^*,v\ra>0\;\mbox{ and }\;\la u^*,v\ra\le 0\;\mbox{ for all }\;
u^*\in\cl\co\bigcup_{i=1}^\infty N(\ox;\O_i).
\end{equation}
Hence for each $\in\N$ we get $\la u^*,v\ra\le 0$ whenever $u^*\in
N(\ox;\O_i)$, which implies that $v\in N^*(\ox;\O_i)$ and therefore $v\in T(\ox;\O_i)$ by \eqref{du1}.
This gives us
$\displaystyle v\in\bigcap_{i=1}^\infty T(\ox;\O_i)$, and so
$\displaystyle x^*\not\in\Big(\bigcap_{i=1}^\infty T(\ox;\O_i)\Big)^*$ due to $\la x^*,v\ra>0$ in \eqref{v}.
Thus we get the inclusion ``$\subset$" in \eqref{eq:CHIPchar1}, which holds as equality. Similar arguments
justify the fulfillment of the parallel duality relationship
\begin{equation}\label{eq:CHIPchar2}
\bigcap_{i=1}^\infty T(\ox;\O_i)=\Big(\cl\co\bigcup_{i=1}^\infty N(\ox;\O_i)\Big)^*.
\end{equation}
Assuming now that the CHIP in (a) holds and employing \eqref{du} and (\ref{eq:CHIPchar1}) for the set
intersection $\O:=\displaystyle\bigcap_{i=1}^\infty\O_i$, we arrive at the equalities
\begin{equation*}
N(\ox;\O)=T^*(\ox;\O)=\Big(\bigcap_{i=1}^\infty
T(\ox;\O_i)\Big)^*=\cl\co\bigcup_{i=1}^\infty N(\ox;\O_i),
\end{equation*}
which give the asymptotic strong CHIP in (b). Conversely, assume that (b) holds. Then employing \eqref{du1} and
(\ref{eq:CHIPchar2}) implies the relationships
\begin{equation*}
T(\ox;\O)=N^*(\ox;\O)=\Big(\cl\co\bigcup_{i=1}^\infty
N(\ox;\O_i)\Big)^*=\bigcap_{i=1}^\infty T(\ox;\O_i),
\end{equation*}
which ensure the fulfillment of the CHIP in (a) and thus establish the equivalence the properties in (a) and (b).
Since the strong CHIP implies the asymptotic strong CHIP due to the closedness of $N(\ox;\O)$, it also implies the
CHIP. The converse implication does not hold even for finitely many sets; counterexamples are presented, in particular,
in \cite{Bauschke-Borwein-Li99,Ernst-Thera07}. $\h$\vspace*{0.05in}

The following simple consequence of Theorem~\ref{Thm:CHIPchar} computes the normal cone to set of feasible solutions
in linear semi-infinite programming with countable inequality constraints; cf.\ \cite{CLMParra09a}.

\begin{Corollary}\label{linear-SIP} {\bf (normal cone to sets of feasible solutions of linear semi-infinite programs
with countable constraints).} Consider the set
\begin{equation}\label{lin}
\O:=\big\{x\in\R^n\big|\;\la a_i,x\ra\le 0,\;i\in\N\big\},
\end{equation}
where the vectors $a_i\in\R^n$ are fixed. Then the normal cone to $\O$ at the origin is computed by
\begin{equation}\label{lin1}
N(0;\O)=\cl\co\Big[\bigcup_{i=1}^\infty\big\{\lm a_i\big|\;\lm\ge 0\}\Big].
\end{equation}
\end{Corollary}
{\bf Proof.} It is easy to see that the set \eqref{lin} is
represented as a countable intersection of sets having the CHIP.
Furthermore, the asymptotic strong CHIP for this system is
obviously \eqref{lin1}. Thus the result follows immediately from
Theorem~\ref{Thm:CHIPchar}. $\h$\vspace*{0.05in}

Let us show that the CHIP may be violated in rather simple situations involving finite and infinite intersections of
convex sets defined by inequalities with convex (while nonlinear) functions.

\begin{Example}\label{Ex:CHIPfail1} {\bf (failure of CHIP for finite and infinite intersections of convex sets).}

{\rm {\bf (i)}  First consider the two convex sets
\begin{equation*}
\O_1:=\big\{(x_1,x_2)\in\R^2\big|\;x_2\ge
x_1^2\big\}\;\mbox{ and }\;\O_2:=\big\{(x_1,x_2)\in\R^2\big|\;x_2\le-x_1^2\big\}
\end{equation*}
and their intersection at $\ox=(0,0)$. We have
\begin{equation*}
\O_1\cap\O_2=\{\ox\},\;\;T(\ox;\O_1)=\R\times\R_+,\mbox{ and }\;
T(\ox;\O_2)=\R\times\R_-.
\end{equation*}
Thus the CHIP does not hold in this case, since
\begin{equation*}
T(\ox;\O_1\cap\O_2)=\{(0,0)\}\ne T(\ox;\O_1)\cap T(\ox;\O_2)=\R\times\{0\}.
\end{equation*}

{\bf (ii)} In the next case we have the CHIP violation for the
countable intersection of convex sets, with the intersection set
having nonempty interior. For each $i\in\N$, define $\ph_i(x):=i
x^2$ if $x<0$ and $\ph_i(x):=0$ if $x\ge 0$. Let $\O_i:=\epi\ph_i$
and $\ox=(0,0)$. It is easy to see that
\begin{equation*}
\bigcap_{i=1}^\infty\O_i=\R_+\times\R_+\;\mbox{ and }\;
T(\ox,\O_i)=\R\times\R_+\;\mbox{ for }\;i\in\N.
\end{equation*}
It gives therefore the relationships
\begin{equation*}
T\Big(\ox,\bigcap_{i=1}^\infty\O_i\Big)=\R_+\times\R_+\ne\bigcap_{i=1}^\infty T(\ox;\O_i)=\R\times\R_+,\quad i\in\N,
\end{equation*}
which show that the CHIP fails for this system of sets at the origin.}
\end{Example}

Of course, we cannot expect to extend the equivalence of
Theorem~\ref{Thm:CHIPchar} to intersections of nonconvex sets. In
what follows we are mainly interested in obtaining calculus rules
for generalized normals (i.e., to get results of the ``strong
CHIP" type) using the nonconvex CHIP from Definition~\ref{chip}
(i.e., a calculus rule for tangents) as an appropriate assumption
together with additional qualification conditions. Observe that
the implication CHIP $\Longrightarrow\,$ strong CHIP does hold
even for finite intersections of convex sets; see
Theorem~\ref{Thm:CHIPchar}. \vspace*{0.05in}

To implement this strategy, we first intend to obtain some
sufficient conditions for the CHIP of countable intersections of
nonconvex sets. Note that a number of sufficient conditions for
the CHIP has been proposed for finite intersections of convex
sets, where convex interpolation techniques play a particularly
important role; see
\cite{Bauschke-Borwein-Li99,Chui-Deutsch-Ward90,Deutsch-Li-Ward97,Li-Ng-Pong07}
and the references therein. However, such techniques do not seem
to be useful in nonconvex settings. To proceed in deriving
sufficient conditions for the CHIP of countable nonconvex
intersections, we explore some other
possibilities.\vspace*{0.05in}

Let us start with extending the concept and techniques of linear
regularity in the direction of
\cite{Bauschke-Borwein-Li99,Li-Ng-Pong07,Song-Zang06} to the case
of infinite nonconvex systems; cf.\ various results and
discussions therein on particular cases of linear regularity and
its applications. Given a countable system of closed sets
$\{\O_i\}_{i\in\N}$, we say that it is {\em linearly regular} at
$\ox\in\O:=\bigcap_{i=1}^\infty\O_i$ if there exist a neighborhood
$U$ of $\ox$ and a positive number $C>0$ such that
\begin{equation}\label{lin-reg1}
\dist(x;\O)\le C\sup_{i\in\N}\big\{\dist(x;\O_i)\big\}\;\mbox{ for all }\;x\in U.
\end{equation}

In the next proposition we denote for convenience the distance
function dist$(x;\O)$ by $d_\O(x)$ and employ the standard notion
of equi-convergence for families of functions.

\begin{Proposition}\label{Prop:CHIPsuf1}{\bf (sufficient conditions for CHIP in terms of linear regularity).}
Let $\{\O_i\}_{i\in\N}$ be a countable system of closed sets in
$\R^n$ with the intersection $\O:=\bigcap_{i=1}^\infty\O_i$, and
let $\ox\in\O$. Assume that the system of sets $\{\O_i\}_{i\in\N}$
is linearly regular at $\ox$ with some $C>0$ in \eqref{lin-reg1}
and that the family of functions $\{d_{\O_i}(\cdot)\}_{i\in\N}$ is
equi-directionally differentiable at $\ox$ in the sense that for
any $h\in\R^n$ the functions
\begin{equation*}
\left\{\frac{d_{\O_i}(\ox+th)}{t},\;i\in\N\right\}
\end{equation*}
equi-converge as $t\dn 0$ to the corresponding directional
derivatives $d'_{\O_i}(\ox;h)$ uniformly in $i\in\N$. Then for all
$h\in\R^n$ and the positive constant $C$ from \eqref{lin-reg1} we
have the estimate
\begin{equation}\label{lin-reg2}
\dist(h;\Lm)\le C\sup_{i\in\N}\big\{\dist(h;\Lm_i)\big\}\;\mbox{
with }\;\Lm:=T(\ox;\O)\;\mbox{ and }\;\Lm_i:=T(\ox;\O_i)\;\mbox{
as }\;i\in\N.
\end{equation}
In particular, the set system $\{\O_i\}_{i\in\N}$ satisfies the CHIP at $\ox$.
\end{Proposition}
{\bf Proof.} Fixing $h\in\R^n$ and using definition \eqref{tc} of the tangent cone, we get
\begin{equation*}
\dist(h;\Lm)=\liminf_{t\dn 0}\dist\Big(h;\frac{\O-\ox}{t}\Big)=
\liminf_{t\dn 0}\frac{\dist(\ox+th;\O)}{t}.
\end{equation*}
When $t$ is small, by the assumed linear regularity yields that
\begin{equation*}
\frac{\dist(\ox+th;\O)}{t}\le C\sup_{i\in\N}\frac{\dist(\ox+th;\O_i)}{t}.
\end{equation*}
Applying \cite[Theorem~4]{BurkeFQ92} with the assumption of equi-directional differentiability, we have
\begin{equation*}
\frac{\dist(\ox+th;\O_i)}{t}\to d'_{\O_i}(\ox;h)=\dist(h;\Lm_i)\;\mbox{ uniformly in }\;i\;\mbox{ as }\;t\dn 0,
\end{equation*}
i.e., for any $\ve>0$ there exists $\delta>0$ such that whenever $t\in(0,\delta)$ we have
\begin{equation*}
\Big\|\frac{\dist(\ox+th;\O_i)}{t}-\dist(h;\Lm_i)\Big\|\le\ve\;\mbox{ for all }\;i\in\N.
\end{equation*}
Hence it holds for any $t\in(0,\delta)$ that
\begin{equation*}
\sup_{i\in\N}\frac{\dist(\ox+th;\O_i)}{t}\le\sup_{i\in\N}\big\{\dist(h;\Lm_i)\big\}+\ve.
\end{equation*}
Combining all the above, we get the estimates
\begin{equation*}
\dist(h;\Lm)\le C\liminf_{t\dn 0}\sup_{i\in\N}\frac{\dist(\ox+th;\O_i)}{t}\le C\sup_{i\in\N}\big\{\dist(h;\Lm_i)\big\}
+C\ve,
\end{equation*}
which imply \eqref{lin-reg2}, since $\ve$ was chosen arbitrarily.
Finally, the CHIP of the system $\{\O_i\}_{i\in\N}$ at $\ox$
follows directly from \eqref{lin-reg2} and the definitions.
$\h$\vspace*{0.05in}

Now we present a consequence of Proposition~\ref{Prop:CHIPsuf1}
that simplifies the verification of linear regularity for
countable set systems.

\begin{Corollary}\label{simpl} {\bf (CHIP via simplified linear regularity).}
Let $\{\O_i\}_{i\in\N}$ be a countable  system of closed subsets in $\R^n$, and let
$\ox\in\O=\bigcap_{i=1}^\infty\O_i$. Assume that the family $\{d(\cdot;\O_i)\}_{i\in\N}$ is
equi-directionally differentiable at $\ox$ and that there are numbers $C>0$, $j\in\N$,
and a neighborhood $U$ of $\ox$ such that we have the estimate
\begin{equation*}
\dist(x;\O)\le C\sup_{i\ne j}\big\{\dist(x;\O_i)\big\}\;\mbox{ for all x }\;\in\O_j\cap U.
\end{equation*}
Then the set system $\{\O_i\}_{i\in\N}$ satisfies the CHIP at $\ox$.
\end{Corollary}
{\bf Proof.} Employing Proposition~\ref{Prop:CHIPsuf1}, it
suffices to show that the set system $\{\O_i\}_{i\in\N}$ is
linearly regular at $\ox$. To proceed, take $r>0$ so small that
\begin{equation*}
\dist(x;\O)\le C\sup_{i\ne j}\big\{\dist(x;\O_i)\big\}\;\mbox{ for all }\;x\in\O_j\cap(\ox+3r\B).
\end{equation*}
Since the distance function is nonexpansive, for every $y\in\O_j\cap(\ox+3r\b)$ and $x\in\R^n$ we have
\begin{equation*}
\begin{aligned}
0&\le C\sup_{i\ne j}\big\{\dist(y;\O_i)\big\}-\dist(y;\O)\le C\sup_{i\ne j}\Big(\big\{\dist(x;\O_i)\big\}+\|x-y\|\Big)-\dist(x;\O)+\|x-y\|\\
&\le C\sup_{i\ne j}\big\{\dist(x;\O_i)\big\}-\dist(x;\O)+(C+1)\|x-y\|.
\end{aligned}
\end{equation*}
Then it follows for all $x\in\R^n$ that
\begin{equation*}
\dist(x;\O)\le(2C+1)\max\big[\sup_{i\ne j}\big\{\dist(x;\O_i)\big\},\;\dist\big(x;\O_j\cap(\ox+3r\B)\big)\big].
\end{equation*}
Thus the linear regularity of $\{\O_i\}_{i\in\N}$ at $\ox$ in the form of
\begin{equation*}
\dist(x;\O)\le(2C+1)\sup_{i\in\N}\big\{\dist(x;\O_i)\big\}
\end{equation*}
would follow now from the relationship
\begin{equation}\label{di}
\dist\big(x;\O_j\cap(\ox+3r\B)\big)=\dist(x;\O_j)\;\mbox{ for all }\;x\in(\ox+r\B).
\end{equation}
To show \eqref{di}, fix a vector $x\in(\ox+r\B)$ above and pick any
$y\in\O_j\setminus(\ox+3r\B)$. This readily gives us $\|x-y\|\ge\|y-\ox\|-\|\ox-x\|\ge 3r-r=2r$ and implies that

\begin{equation*}
\dist\big(x;\O_j\setminus(\ox+3r\B)\big)\ge 2r\;\mbox{ while }\;\dist\big(x;\O_j\cap(\ox+3r\B)\big)\le\|x-\ox\|\le r.
\end{equation*}
Hence we get the equalities
\begin{equation*}
\dist(x;\O_j)=\min\big\{\dist\big(x;\O_j\setminus(\ox+3r\B)\big),\;\dist\big(x;\O_j\cap(\ox+3r\B)\big)\big\}=
\dist\big(x;\O_j\cap (\ox+3r\B)\big),
\end{equation*}
which justify \eqref{di} and thus complete the proof of the corollary. $\h$\vspace*{0.05in}

The next proposition, which holds in fact for arbitrary (not only
countable) intersections of sets, establishes a new sufficient
condition for the CHIP of $\{\O_i\}_{i\in\N}$. To formulate it, we
introduce a notion of the {\em tangential rank} of the
intersection $\O:=\bigcap_{i=1}^\infty\O_i$ at $\ox\in\O$ by
\begin{equation}\label{rank}
\rho_\O(\ox):=\inf_{i\in\N}\left\{\limsup_{x\to\ox\atop
x\in\O_i\setminus\{\ox\}}\frac{\dist(x;\O)}{\|x-\ox\|}\right\},
\end{equation}
where we put $\rho_\O(\ox):=0$ if $\O_i=\{\ox\}$ for at least one $i\in\N$.

\begin{Proposition}\label{Prop:CHIPsuf2} {\bf (sufficient condition for CHIP via tangential rank of intersection).}
Given a countable system of  closed sets $\{\O_i\}_{i\in\N}$ in
$\R^n$, suppose that $\rho_\O(\ox)=0$ for the tangential rank of
their intersection $\O:=\bigcap_{i=1}^\infty\O_i$ at $\ox\in\O$.
Then this system exhibits the CHIP at $\ox$.
\end{Proposition}
{\bf Proof}. The result holds trivially if $\O_i=\{\ox\}$ for some $i\in\N$. Assume that
$\O_i\setminus\{\ox\}\ne\emp$ for all $i\in\N$ and observe that
$T(\ox;\O)\subset T(\ox;\O_i)$ whenever $i\in\N$. Thus we have
\begin{equation*}
T(\ox;\O)\subset\bigcap_{i\in\N}T(\ox;\O_i).
\end{equation*}
To prove the reverse inclusion, fix an arbitrary vector $0\ne v\in\bigcap_{i=1}^\infty T(\ox;\O_i)$. By the assumption
of $\rho_\O(\ox)=0$ and definition \eqref{rank}, for any $k\in\N$ we find a set $\O_k$ from the system under
consideration such that
\begin{equation*}
\limsup_{x\to\ox\atop x\in\O_k\setminus\{\ox\}}\frac{\dist(x;\O)}{\|x-\ox\|}<\frac{1}{k}.
\end{equation*}
Since $v\in T(\ox;\O_k)$, there are sequences $\{x_j\}_{j\in\N}\subset\O_k$ and $t_j\dn 0$ satisfying
\begin{equation*}
x_j\to\ox\quad{\rm and}\quad\frac{x_j-\ox}{t_j}\to v\quad{\rm as}\;j\to\infty,
\end{equation*}
which in turn implies the limiting estimate
\begin{equation*}
\limsup_{j\to\infty}\frac{\dist(x_j;\O)}{\|x_j-\ox\|}<\frac{1}{k}.
\end{equation*}
The latter allows us to find a vector $x_k\subset\{x_j\}_{j\in\N}$
with $\|x_k-\ox\|\le 1/k$ and the corresponding number $t_k\le
1/k$ such that
\begin{equation*}
\left\|\frac{x_k-\ox}{t_k}-v\right\|\le\frac{1}{k}\quad{\rm and}\quad
\frac{\dist(x_k;\O)}{\|x_k-\ox\|}<\frac{1}{k}.
\end{equation*}
Then it follows that there exists $z_k\in\O$ satisfying the relationships
\begin{equation*}
\|z_k-x_k\|<\frac{1}{k}\|x_k-\ox\|\le\frac{1}{k^2}.
\end{equation*}
Combining all the above together gives us the estimates
\begin{equation*}
\left\|\frac{z_k-\ox}{t_k}-v\right\|\le
\left\|\frac{z_k-x_k}{t_k}\right\|+\left\|\frac{x_k-\ox}{t_k}-v\right\|
\le\frac{1}{k}\left\|\frac{x_k-\ox}{t_k}\right\|+\frac{1}{k}
\le\frac{1}{k}\Big(\|v\|+\frac{1}{k}\Big)+\frac{1}{k},\quad k\in\N.
\end{equation*}
Now letting $k\to\infty$, we get $z_k\st{\O}{\longrightarrow}\ox$,
$t_k\dn0$, and
$\disp\left\|\frac{z_k-\ox}{t_k}-v\right\|\longrightarrow 0$. The
latter verifies that $v\in T(\ox;\O)$ and thus completes the proof
of the proposition. $\h$\vspace*{0.05in}

To conclude our discussions on the CHIP, we give yet another
verifiable condition ensuring the fulfillment of this property for
countable intersections of closed sets. We say that a set $A$ is
of {\em invex type} if it can be represented as the complement to
a union with respect to $t\in T$ of some open convex sets $A_t$,
i.e.,
\begin{equation}\label{inv0}
A=\R^n\setminus\bigcup_{t\in T}A_t,
\end{equation}
The following lemma needed for the next proposition is also used in Section~5.

\begin{Lemma}\label{Lem:Invex} {\bf (sets of invex type).}  Let $A\subset\R^n$ be a set of invex type, and
let $\ox\in\bigcap_{t\in T}\bd A_t\cap\bd A$ be taken from the boundary intersections. Then we have the inclusion involving the tangent cone $T(\ox;A)$ to $A$ at $\ox$:
\begin{equation}\label{inv}\ox+T(\ox;A)\subset A.
\end{equation}
\end{Lemma}
{\bf Proof}. To justify inclusion \eqref{inv}, suppose on the contrary that there is $v\in T(\ox;A)$ such that
$\ox+v\notin A$. For this vector $v$ we find by definition \eqref{tc} sequences $s_k\dn 0$ and $x_k\in A$ such that
$\frac{x_k-\ox}{s_k}\to v$ as $k\to\infty$. Since $\ox+v\notin A$, by invexity \eqref{inv0} there exists an index
$t\in T$ for which $\ox+v\in A_t$. Thus we get the inclusion
\begin{equation*}
\ox+\frac{x_k-\ox}{s_k}\in A_t\;\mbox{ for all }\;k\in\N\;\mbox{ sufficiently large}.
\end{equation*}
Then employing the convexity of $A_t$ gives us that
\begin{equation*}
x_k=(1-s_k)\ox+s_k\Big(\ox+\frac{x_k-\ox}{s_k}\Big)\in A_t
\end{equation*}
for the fixed index $t\in T$ and all large numbers $k\in\N$. This contradicts the choice of $x_k\in A$ and thus
justifies the claimed inclusion \eqref{inv}. $\h$\vspace*{0.05in}.

Now we are ready to derive the aforementioned sufficient condition for the CHIP.

\begin{Proposition}\label{Prop:CHIPsuf3} {\bf (CHIP for countable intersections of invex-type sets).} Given a
countable system $\{\O_i\}_{i\in\N}$ in $\R^n$, assume that there is a $($possibly infinite$)$ index subset
$J\subset\N$ such that each $\O_i$ for $i\in J$ is the complement to an open and convex set in $\R^n$ and that
\begin{equation}\label{inv2}
\ox\in\bigcap_{i\in J}\bd\O_i\cap\inter\bigcap_{i\not\in J}\O_i
\end{equation}
for some $\ox$. Then the system $\{\O_i\}_{i\in\N}$ enjoys the CHIP at $\ox$.
\end{Proposition}
{\bf Proof}. Take any $\O_i$ with $i\in J$ and find a convex and open set $A\subset\R^n$ such that $\O=R^n\setminus A$.
Then $\ox\in\bd A\cap\bd\O_i$ by \eqref{inv2}. Then Lemma~\ref{Lem:Invex} ensures that $\ox+T(\ox;\O_i)\subset\O_i$ for
this index $i\in J$. By the choice of $\ox$ in \eqref{inv2} we have furthermore that
\begin{equation*}
\bigcap_{i=1}^\infty T(\ox;\O_i)=\bigcap_{i\in J}T(\ox;\O_i)\subset\bigcap_{i\in J}(\O_i-\ox).
\end{equation*}
Since the set on the left-hand side of the latter inclusion is a cone, it follows that
\begin{equation}\label{inv3}
\bigcap_{i=1}^\infty T(\ox;\O_i)\subset T\Big(0;\bigcap_{i\in J}(\O_i-\ox)\Big)=
T\Big(\ox;\bigcap_{i\in J}\O_i\Big)=T\Big(\ox;\bigcap_{i=1}^\infty\O_i\Big).
\end{equation}
As the opposite inclusion in \eqref{inv3} is obvious, we conclude that the CHIP is satisfied for the countable set
system $\{\O_i\}_{i\in\N}$ at $\ox$. $\h$\vspace*{0.1in}

In the last part of this section we show that the CHIP for countable intersections of nonconvex sets, combined with
some other classification conditions, allows us to derive principal calculus rules for representing
{\em generalized normals to infinite set intersections}. Thus the verifiable sufficient conditions for the CHIP
established above largely contribute to the implementation of these calculus rules. Note that the results obtained in
this direction provide new information even for convex set intersections, since in this case they furnish the required
implication CHIP $\,\Longrightarrow\,$ strong CHIP, which does not hold in general nonconvex settings; see
Theorem~\ref{Thm:CHIPchar} for more discussions.\vspace*{0.05in}

First we formulate and discuss appropriate qualification conditions for countable systems of sets in terms of the basic
normal cone \eqref{nc}.

\begin{Definition}\label{qc} {\bf (normal closedness and qualification conditions for countable set systems).}
Let $\{\O_i\}_{i\in\N}$ be a countable system of sets, and let $\ox\in\bigcap_{i=1}^\infty\O_i$. We say that:

{\bf (a)} The set system $\{\O_i\}_{i\in\N}$ satisfies the {\sc normal closedness condition} {\rm (NCC)} at $\ox$ if
the combination of basic normals
\begin{equation}\label{eq:CC}
\Big\{\sum_{i\in I}x^*_i\Big|\;x^*_i\in N(\ox;\O_i),\;I\in{\cal L}\Big\}\;\mbox{ is closed in }\;\R^n,
\end{equation}
where ${\cal L}$ stands for the collection of all the finite subsets of $\N$.

{\bf (b)} The system $\{\O_i\}_{i\in\N}$ satisfies the {\sc normal qualification condition} {\rm (NQC)} at $\ox$
if the following implication holds:
\begin{equation}\label{eq:QC-N(O)}
\left[\sum_{i=1}^{\infty}x^*_i=0,\;\; x^*_i\in
N(\ox;\O_i)\right]\Longrightarrow\Big[x^*_i=0\;\mbox{ for all }\;
i\in\N\Big].
\end{equation}
\end{Definition}

The NCC in Definition~\ref{qc}(a) relates to various versions of
the so-called {\em Farkas-Minkowski qualification condition} and
its extensions for finite and infinite systems of sets. We refer
the reader to, e.g., \cite{DMN09b,DMN09a} and the bibliographies
therein, as well as to subsequent discussions in Section~4, for a
number of results in this direction concerning convex infinite
inequality systems and to \cite{CLMParra09b} for more details on
linear inequality systems with arbitrary index sets in general
Banach spaces.

The NQC in Definition~\ref{qc}(b) is a direct extension of the corresponding condition \eqref{eq:QC-N(Lm)}) for system
of cones. The counterpart of \eqref{eq:QC-N(O)} for finite systems of sets is studied and applied in
\cite{m-book1,m-book2} under the same name. The following proposition presents a simple sufficient condition for the
validity of the NQC in the case of countable systems of convex sets.

\begin{Proposition}\label{Prop:QCcvex}{\bf (NQC for countable systems of convex sets).} Let $\{\O_i\}_{i\in\N}$ be a
system of convex sets for which there is an index $i_0\in\N$ such that
\begin{equation}\label{eq:QCcvex}
\O_{i_0}\cap\bigcap_{i\ne i_0}\inter\O_i\ne\emp.
\end{equation}
Then the NQC in \eqref{eq:QC-N(O)} is satisfied for the system $\{\O_i\}_{i\in\N}$ at any $\ox\in\bigcap_{i=1}^
\infty\O_i$.
\end{Proposition}
{\bf Proof.} Suppose without loss of generality that $i_0=1$ and fix some $w\in\O_1\cap\bigcap_{i=2}^\infty\inter\O_i$.
Taking any normals $x^*_i\in N(\ox;\O_i)$ as $i\in\N$ satisfying
\begin{equation*}
\sum_{i=1}^{\infty}x^*_i=0,
\end{equation*}
we get by the convexity of the sets $\O_i$ that $\la x^*_i,w-\ox\ra\le 0$ for all $i\in\N$. Then it follows that
\begin{equation*}
\la x^*_i,w-\ox\ra=-\sum_{j\ne i}\la x^*_j,w-\ox\ra\ge 0,\quad i\in\N,
\end{equation*}
which yields $\la x^*_i,w-\ox\ra=0$ whenever $i\in\N$. Next fix $\ve>0$ and find $m\in\N$ so large that
\begin{equation*}
\left\|\sum_{i=m+1}^{\infty}x^*_i\right\|\le\ve.
\end{equation*}
Pick $u\in\R^n$ with $\|u\|=1$ and taking into account that
$w\in\bigcap_{i=2}^m\inter\O_i$, we get
\begin{equation*}
\lm\la x^*_i,u\ra=\la x^*_i,w+\lm u-\ox\ra\le 0,\quad i=2,3,\ldots,
\end{equation*}
whenever $\lm>0$ is sufficiently small. This implies that
\begin{equation*}
\lm\la x^*_1,u\ra=-\lm\sum_{i=2}^m\la x^*_i,u\ra-\lm\sum_{i=m+1}^\infty\la x^*_i,u\ra
\ge-\lm\left\|\sum_{i=m+1}^{\infty}x^*_i\right\|\cdot\|u\|\ge-\lm\ve,
\end{equation*}
which gives $\la x^*_1,u\ra\ge-\ve$. Since $\ve>0$ was chosen arbitrarily, we conclude that $\la x^*_1,u\ra\ge 0$.
Repeating the same procedure for $-u$ shows that $\la x^*_1,-u\ra\ge 0$ and so
$\la x^*_1,u\ra=0$ for all $u\in\R^n$ with $\|u\|=1$. This implies that $x^*_1=0$.
The same procedure ensures that $x^*_i=0$ for all $i\in\N$, which completes the proof
of the proposition. $\h$\vspace*{0.05in}

Finally, we obtain the main result of this section, which
expresses Fr\'echet normal to infinite set intersections via basic
normals to the sets involved under the above CHIP and
qualification conditions. This major calculus rule for arbitrary
closed sets employs the corresponding intersection rule for cones
from Theorem~\ref{Thm:FnorLm}, which is based on the tangential
extremal principle.

\begin{Theorem}\label{Thm:FnorO}{\bf (generalized normals to countable set intersections).} Let $\{\O_i\}_{i\in\N}$
be a countable system of closed sets in $\R^n$, and let $\ox\in\O:=\bigcap_{i=1}^\infty\O_i$. Assume that the CHIP in
\eqref{chip} and NQC in \eqref{eq:QC-N(O)} are satisfied for $\{\O_i\}_{i\in\N}$ at $\ox$.
Then we have the inclusion
\begin{equation}\label{f-inter}
\Hat N(\ox;\O)\subset\cl\Big\{\sum_{i\in I}x^*_i\Big|\;x^*_i\in N(\ox;\O_i),\;I\in{\cal L}\Big\},
\end{equation}
where ${\cal L}$ stands for the collection of all the finite subsets of $\N$. If in addition the CQC in \eqref{eq:CC}
holds for $\{\O_i\}_{i\in\N}$ at $\ox$, then the closure operation can be omitted on the right-hand side of
\eqref{f-inter}.
\end{Theorem}
{\bf Proof}. Using the assumed CHIP for $\{\O_i\}_{i\in\N}$ at $\ox$, constructions \eqref{tc} and \eqref{eq:e-nor},
and the duality correspondence \eqref{du} gives us
\begin{equation}\label{f-int}
\Hat N(\ox;\O)=\Hat N\big(0;T(\ox;\O)\big)=\Hat
N\Big(0;\bigcap_{i=1}^\infty T(\ox;\O_i)\Big).
\end{equation}
It follows from \eqref{eq:N(Lm)-N(O)} that $N\big(0;T(\ox;\O_i)\big)\subset N(\ox;\O_i)$ for all $i\in\N$, and thus
the assumed NQC in \eqref{eq:QC-N(O)} implies the conic one in (\ref{eq:QC-N(Lm)}). Applying Theorem~\ref{Thm:FnorLm},
we have
\begin{equation*}
\Hat N\Big(0;\bigcap_{i=1}^\infty T(\ox;\O_i)\Big)\subset\cl\Big\{\sum_{i\in I}x^*_i\Big|\;x^*_i\in N\big(0;T(\ox;\O_i)
\big),\;I\in{\cal L}\Big\}.
\end{equation*}
Now the intersection rule \eqref{f-inter} follows from \eqref{eq:N(Lm)-N(O)} and \eqref{f-int}. Finally, the
closure operation in \eqref{f-inter} can be obviously dropped if the system $\{\O_i\}_{i\in\N}$ satisfies the CQC at
$\ox$. $\h$

\section{Applications to Semi-Infinite Programming}
\setcounter{equation}{0}

This section is devoted to deriving necessary optimality
conditions for various problems of semi-infinite programming (SIP)
with {\em countable} constraints. As mentioned in Section~1,
problems with countable constraints are among the most difficult
in SIP, in comparison with conventional ones involving constraints
indexed by compact sets. In fact, SIP problems with countable
constraints are not different from seemingly more general problems
with arbitrary index sets. Problems of the latter class have drawn
particular attention in a number of recent publications, where
some special structures of this type (mostly with linear and
convex inequality constraints) have been considered; see, e.g.,
\cite{CLMParra09b,DMN09b,DMN09a} and the references therein. In
this section we derive, based on the tangential extremal principle
and its calculus consequences, new optimality conditions for SIP
with various types of countable constraints and compare them with
those known in the literature.\vspace*{0.05in}

Let us start with SIP involving countable constraints of the {\em geometric type}:
\begin{equation}\label{eq:SIP1}
\mbox{minimize }\;\ph(x)\;\mbox{ subject to }\;x\in\O_i\;\mbox{ as }\;i\in\N,
\end{equation}
where $\ph\colon\R^n\to\oR$ is an extended-real-valued function, and where $\{\O_i\}_{i\in\N}\subset\R^n$ is a
countable system of constraint sets.
Considering in general problems with nonsmooth and nonconvex cost functions and following the classification of
\cite[Chapter~5]{m-book2}, we derive necessary optimality conditions of two kinds for \eqref{eq:SIP1} and other SIP
{\em minimization} problems: {\em lower subdifferential} and {\em upper subdifferential} ones. Conditions of the
``lower" kind are more conventional for minimization dealing with usual (lower) subdifferential constructions.
On the other hand, conditions of the ``upper" kind employ upper subdifferential (or superdifferential) constructions,
which seem to be more appropriate for maximization problems while bringing significantly stronger information for
special classes of minimizing cost functions in comparison with lower subdifferential ones; see \cite{m-book2}
for more discussions, examples, and references.

We begin with upper subdifferential optimality conditions for \eqref{eq:SIP1}. Given $\ph\colon\R^n\to\oR$ finite at
$\ox$, the {\em upper subdifferential} of $\ph$ at $\ox$ used in this paper is of the Fr\'echet type defined by
\begin{equation}\label{u-fr}
\Hat\partial^+\ph(\ox):=-\Hat\partial(-\ph)(\ox)=\Big\{x^*\in\R^n\Big|\;\limsup_{x\to\ox}\frac{\ph(x)-\ph(\ox)-\la x^*,
x-\ox\ra}{\|x-\ox\|}\le 0\Big\}
\end{equation}
via \eqref{fr}. Note that $\Hat\partial^+\ph(\ox)$ reduces to the
upper subdifferential (or superdifferential) of convex analysis if
$\ph$ is concave. Furthermore, the subdifferential sets
$\Hat\partial\ph(\ox)$ and $\Hat\partial^+\ph(\ox)$ are nonempty
simultaneously if and only if $\ph$ is Fr\'echet differentiable at
$\ox$.\vspace*{0.05in}

As before, in the next theorem and in what follows the symbol ${\cal L}$ stands for the collection of all the finite
subsets of the natural series $\N$.

\begin{Theorem}\label{Thm:NecPg1}{\bf (upper subdifferential conditions for SIP with countable geometric constraints).}
Let $\ox$ be a local optimal solution to problem \eqref{eq:SIP1},
where $\ph\colon\R^n\to\oR$ is an arbitrary extended-real-valued
function finite at $\ox$, and where the sets $\O_i\subset\R^n$ for
$i\in\N$ are locally closed around $\ox$. Assume that the system
$\{\O_i\}_{i\in\N}$ has the CHIP at $\ox$ and satisfies the NQC of
Definition~{\rm\ref{qc}(b)} at this point. Then we have the set
inclusion
\begin{equation}\label{u-fr1}
-\Hat\partial^+\ph(\ox)\subset\cl\Big\{\sum_{i\in I}x^*_i\Big|\;
x^*_i\in N(\ox;\O_i),\;I\in{\cal L}\Big\},
\end{equation}
which reduces to that of
\begin{equation}\label{u-fr2}
0\in\nabla\ph(\ox)+\cl\Big\{\sum_{i\in I}x^*_i\Big|\;
x^*_i\in N(\ox;\O_i),\;I\in{\cal L}\Big\}.
\end{equation}
if $\ph$ is Fr\'echet differentiable at $\ox$. If in addition the NCC of Definition~{\rm\ref{qc}(a)} holds for
$\{\O_i\}_{i\in\N}$ at $\ox$, then the closure operations can be omitted in \eqref{u-fr1} and \eqref{u-fr2}.
\end{Theorem}
{\bf Proof}. It follows from \cite[Proposition~5.2]{m-book2} that
\begin{equation}\label{u-fr3}
-\Hat\partial^+\ph(\ox)\subset\Hat N\Big(\ox;\bigcap_{i=1}^\infty\O_i\Big).
\end{equation}
Applying now to \eqref{u-fr3} the representation of Fr\'echet
normals to countable set intersections from
Theorem~\ref{Thm:FnorO} under the assumed CHIP and NQC, we arrive
at \eqref{u-fr1}, where the closure operation can be omitted when
the NCC holds at $\ox$. If $\ph$ is Fr\'echet differentiable at
$\ox$, it follows that
$\Hat\partial^+\ph(\ox)=\{\nabla\ph(\ox)\}$, and thus
\eqref{u-fr1} reduces to \eqref{u-fr2}. $\h$\vspace*{0.05in}

Note that the set inclusion \eqref{u-fr1} is trivial if
$\Hat\partial^+\ph(\ox)=\emp$, which is the case of, e.g.,
nonsmooth convex functions. On the other hand, the upper
subdifferential necessary optimality condition \eqref{u-fr1} may
be much more selective than its lower subdifferential counterparts
when $\Hat\partial^+\ph(\ox)\ne\emp$, which happens, in
particular, for some remarkable classes of functions including
concave, upper regular, semiconcave, upper-$C^1$, and other ones
important in various applications. The reader can find more
information and comparison in \cite[Subsection~5.1.1]{m-book2} and
the commentaries therein concerning problems with finitely many
geometric constraints. \vspace*{0.05in}

Next let us present a lower subdifferential condition for the SIP
problem \eqref{eq:SIP1} involving the basic subdifferential
\eqref{sub}, which is nonempty for majority of nonsmooth
functions; in particular, for any local Lipschitzian one. To
formulate this condition, recall the notion of the {\em singular
subdifferential} of $\ph$ at $\ox$ defined by
\begin{equation}\label{sin}
\partial^\infty\ph(\ox):=\big\{x^*\in\R^n\big|\;(x^*,0)\in N\big((\ox;\ph(\ox));\epi\ph\big)\big\}.
\end{equation}
Note that $\partial^\infty\ph(\ox)=\{0\}$ if $\ph$ is locally
Lipschitzian around $\ox$. Recall also that a set $\O$ is {\em
normally regular} at $\ox$ if $N(\ox;\O)=\Hat N(\ox;\O)$. This is
the case, in particular, of locally convex and other ``nice" sets;
see, e.g., \cite{m-book1,Rockafellar-Wets-VA} and the references
therein.

\begin{Theorem}\label{Thm:Nec2} {\bf (lower subdifferential subdifferential conditions for SIP with countable
geometric constraints.)} Let $\ox$ be a local optimal solution to
problem \eqref{eq:SIP1} with a lower semicontinuous cost function
$\ph\colon\R^n\to\oR$ finite at $\ox$ and a countable system
$\{\O_i\}_{i\in\N}$ of sets locally closed around $\ox$. Assume
that the feasible solution set $\O:=\bigcap_{i=1}^\infty\O_i$ is
normally regular at $\ox$, that the system $\{\O_i\}_{i\in\N}$
satisfies the CHIP \eqref{chip} and the NQC \eqref{eq:QC-N(O)} at
$\ox$, and that
\begin{equation}\label{sin1}
\cl\Big\{\sum_{i\in I}x^*_i\Big|\;
x^*_i\in N(\ox;\O_i),\;I\in{\cal L}\Big\}\bigcap\big(-\partial^\infty\ph(\ox)\big)=\{0\},
\end{equation}
which holds, in particular, when $\ph$ is locally Lipschitzian around $\ox$. Then we have
\begin{equation}\label{lower1}
0\in\partial\ph(\ox)+\cl\Big\{\sum_{i\in I}x^*_i\Big|\;
x^*_i\in N(\ox;\O_i),\;I\in{\cal L}\Big\}.
\end{equation}
The closure operations can be omitted in \eqref{sin1} and
\eqref{lower1} if the NCC \eqref{eq:CC} is satisfied at $\ox$.
\end{Theorem}
{\bf Proof}. It follows from \cite[Proposition~5.3]{m-book2} that
\begin{equation}\label{sin2}
0\in\partial\ph(\ox)+N(\ox;\O\;\mbox{ provided that }\;\partial^\infty\ph(\ox)\cap\big(-N(\ox;\O)\big)=\{0\}
\end{equation}
for the optimal solution $\ox$ to the problem under consideration
with the feasible solution set $\O:=\bigcap_{i=1}^\infty\O_i$.
Since the set $\O$ is normally regular at $\ox$, we can replace
$N(\ox;\O)$  by $\Hat N(\ox;\O)$ in \eqref{sin2}. Applying now
Theorem~\ref{Thm:FnorO} to the countable set intersection $\O$ in
\eqref{sin2} under the assumptions made, we arrive at all the
conclusions of this theorem. $\h$\vspace*{0.05in}

Next we consider a SIP problem with {\em countable operator constraints} defined by:
\begin{equation}\label{sip-op}
\mbox{minimize }\;\ph(x)\;\mbox{ subject
to }\;f(x)\in\Th_i\;\mbox{ as }\;i\in\N,
\end{equation}
where $\ph\colon\R^n\to\oR$, $\Th_i\subset\R^m$ for $i\in\N$, and
$f\colon\R^n\to\R^m$.  The following statements are consequences
of Theorems~\ref{Thm:NecPg1} and \ref{Thm:Nec2}, respectively.

\begin{Corollary}\label{sip-op1} {\bf (upper and lower subdifferential conditions for SIP with operator constraints).}
Let $\ox$ be a local optimal solution to \eqref{sip-op}, where the
function $\ph\colon\R^n$ is finite at $\ox$, where the mapping
$f\colon\R^n\to\R^m$ is strictly differentiable at $\ox$ with the
surjective $($full rank$)$ derivative, and where the sets
$\Th_i\subset\R^m$ as $i\in\N$ are locally closed around $f(\ox)$
while satisfying the CHIP \eqref{chip} and NQC \eqref{eq:QC-N(O)}
conditions at this point. The following assertions holds:

{\bf (i)} We have the upper subdifferential optimality condition:
\begin{equation}\label{u-op}
-\Hat\partial^+\ph(\ox)\subset\cl\Big\{\sum_{i\in I}\nabla f(\ox)^*
y^*_i\Big|\;y^*_i\in N\big(f(\ox);\Th_i\big),\;I\in{\cal L}\Big\},
\end{equation}

{\bf (ii)} If $\ph$ is lower semicontinuous around $\ox$ and
\begin{equation}\label{sin3}
\cl\Big\{\sum_{i\in I}\nabla f(\ox)^*
y^*_i\Big|\;y^*_i\in N\big(f(\ox);\Th_i\big),\;I\in{\cal L}\Big\}\bigcap\big(-\partial^\infty\ph(\ox)\big)=\{0\},
\end{equation}
then we have the inclusion
\begin{equation}\label{sin4}
0\in\partial\ph(\ox)+\cl\Big\{\sum_{i\in I}\nabla f(\ox)^*
y^*_i\Big|\;y^*_i\in N\big(f(\ox);\Th_i\big),\;I\in{\cal L}\Big\}.
\end{equation}
Furthermore, the closure operations can be omitted in
\eqref{u-op}--\eqref{sin4} if  the set system $\{\Th_i\}_{i\in\N}$
satisfies the NCC \eqref{eq:CC} at $f(\ox)$.
\end{Corollary}
{\bf Proof}. Observe that problem \eqref{sip-op} can be
equivalently rewritten in the geometric form \eqref{eq:SIP1} with
$\O_i:=f^{-1}(\Th_i)$, $i\in\N$. Then employing the well-known
results on representing the tangent and normal cones in \eqref{tc}
and \eqref{nc} to inverse images of sets under strict
differentiable mappings with surjective derivatives (see, e.g.,
\cite[Theorem~1.17]{m-book1} and
\cite[Exercise~6.7]{Rockafellar-Wets-VA}), we have
\begin{equation}\label{inver}
T\big(\ox;f^{-1}(\Th)\big)=\nabla
f(\ox)^{-1}T\big(f(\ox);\Th\big)\;\mbox{ and
}\;N\big(\ox;f^{-1}(\Th)\big)=\nabla
f(\ox)^*N\big(f(\ox);\Th\big).
\end{equation}
It follows from the surjectivity of $\nabla f(\ox)$ that the CHIP
and NQC for $\{\Th_i\}_{i\in\N}$ at $f(\ox)$ are equivalent,
respectively, to the CHIP and NQC of $\{\O_i\}_{i\in\N}$ at $\ox$;
see \cite[Lemma~1.18]{m-book1}. This implies the equivalence
between the qualification and optimality conditions
\eqref{u-op}--\eqref{sin4} for problem \eqref{sip-op} under the
assumptions made and the corresponding conditions \eqref{u-fr1},
\eqref{sin1}, and \eqref{lower1} for problem \eqref{eq:SIP1}
established in Theorems~\ref{Thm:NecPg1} and \ref{Thm:Nec2}. To
complete the proof of the corollary, it suffices to observe
similarly to \eqref{inver} that the assumed NCC for
$\{\Th_i\}_{i\in\N}$ at $f(\ox)$ is equivalent under the
surjectivity of $\nabla f(\ox)$ to the NCC \eqref{eq:CC} for the
inverse images $\{\O_i\}_{i\in\N}$ at $\ox$. Thus the possibility
to omit the closure operations in the framework of the corollary
follows directly from the corresponding statements of
Theorems~\ref{Thm:NecPg1} and \ref{Thm:Nec2}. $\h$\vspace*{0.05in}

The rest of this section concerns SIP problems with {\em countable inequality constraints}:
\begin{equation}\label{sip-ine}
\mbox{minimize }\;\ph(x)\;\mbox{ subject
to }\;\ph_i(x)\le 0\;\mbox{ as }\;i\in\N,
\end{equation}
where the cost function $\ph$ is as in problems \eqref{eq:SIP1}
and \eqref{sip-op} while the constraints functions
$\ph_i\colon\R^n\to\oR$, $i\in\N$, are lower semicontinuous around
the reference optimal solution. Note that problems with infinite
inequality constraints are considered in the vast majority of
publications on semi-infinite programming, where the main
attention is paid to the case of convex or linear infinite
inequalities; see below some comparison with known results for SIP
of the latter types.

Although our methods are applied to problems \eqref{sip-ine} of
the general inequality type, for simplicity and brevity we focus
here on the case when the constraint functions $\ph_i$, $i\in\N$,
are locally Lipschitzian around the optimal solution. In the
general case we need to involve the singular subdifferential
\eqref{sin} of these functions; see the proofs below. Let us first
introduce subdifferential counterparts of the normal qualification
and closedness conditions from Definition~\ref{qc}.

\begin{Definition}\label{sub-qc} {\bf (subdifferential closedness and qualification conditions for countable
inequality constraints).} Consider a countable constraint system $\{\O_i\}_{i\in\N}\subset\R^n$ with
\begin{equation}\label{ine}
\O_i:=\big\{x\in\R^n\big|\;\ph_i(x)\le 0\big\},\quad i\in\N,
\end{equation}
where the functions $\ph_i$ are locally Lipschitzian around $\ox\in\bigcap_{i=1}^\infty\O_i$. We say that:

{\bf (a)} The system $\{\O_i\}_{i\in\N}$ in \eqref{ine} satisfies
the {\sc subdifferential closedness condition} $($SCC$)$ at $\ox$
if the set
\begin{equation}\label{scc}
\Big\{\sum_{i\in I}\lm_i\partial\ph_i(\ox)\Big|\;\lm_i\ge 0,\;\lm_i\ph_i(\ox)=0,\;I\in{\cal L}\Big\}\;\mbox{ is
closed in }\;\R^n.
\end{equation}

{\bf (b)} The system $\{\O_i\}_{i\in\N}$ in \eqref{ine} satisfies the {\sc subdifferential qualification condition}
$($SQC$)$ at $\ox$ if the following implication holds:
\begin{equation}\label{eq:QCgi}
\Big[\sum_{i=1}^\infty\lm_i x^*_i=0,\;x^*_i\in\partial\ph_i(\ox),\;\lm_i\ge 0,\;\lm_i\ph_i(\ox)=0\Big|\Longrightarrow
\big[\lm_i=0\;\mbox{ for all }\;i\in\N\big].
\end{equation}
\end{Definition}

The next theorem provides necessary optimality conditions of both upper and lower subdifferential types for SIP problems
\eqref{sip-ine} without any smoothness and/or convexity assumptions.

\begin{Theorem}\label{Thm:NecPf} {\bf (upper and lower subdifferential conditions for general SIP with inequality
constraints).} Let $\ox$ be a local optimal solution to problem \eqref{sip-ine}, where the constraint functions
$\ph_i\colon\R^n\to\oR$ are locally Lipschitzian around $\ox$ for all $i\in\N$. Assume that the level set system
$\{\O_i\}_{i\in\N}$ in \eqref{ine} has the CHIP at $\ox$ and that the SQC \eqref{eq:QCgi} is satisfied at this point.
Then the following assertions hold:

{\bf (i)} We have the upper subdifferential optimality condition:
\begin{equation}\label{up-ine}
-\Hat\partial^+\ph(\ox)\subset\cl\Big\{\sum_{i\in I}\lm_i\partial\ph_i(\ox)\Big|\;\lm_i\ge 0,\;\lm_i\ph_i(\ox)=0,
\;I\in{\cal L}\Big\},
\end{equation}
where the closure operation can be omitted if the SCC \eqref{scc} is satisfied at $\ox$.

{\bf (ii)} Assume in addition that $\ph$ is lower semicontinuous around $\ox$ and that
\begin{equation}\label{lo-qc}
\cl\Big\{\sum_{i\in I}\lm_i\partial\ph_i(\ox)\Big|\;\lm_i\ge 0,\;\lm_i\ph_i(\ox)=0,\;I\in{\cal L}\Big\}
\bigcap\big(-\partial^\infty\ph(\ox)\big)=\{0\},
\end{equation}
which is automatic if $\ph$ is locally Lipschitzian around $\ox$. Then
\begin{equation}\label{lo-sub}
0\in\partial\ph(\ox)+\cl\Big\{\sum_{i\in I}\lm_i\partial\ph_i(\ox)\Big|\;\lm_i\ge 0,\;\lm_i\ph_i(\ox)=0,\;I\in{\cal L}
\Big\}
\end{equation}
with removing the closure operation in \eqref{lo-qc} and \eqref{lo-sub} when the SCC \eqref{scc} holds at $\ox$.
\end{Theorem}
{\bf Proof}. It is well known from the calculus of basic normals and subgradients that
\begin{equation}\label{eq:levelset}
N(\ox;\O)\subset\R_+\partial\vartheta(\ox):=\big\{\lm
x^*\in\R^n\big|\;x^*\in\partial\vartheta(\ox),\;\lm\ge 0\big\}\;
\mbox{ for }\;\O:=\big\{x\in\R^n\big|\;\vartheta(x)\le 0\big\}
\end{equation}
provided that $\vartheta\colon\R^n\to\oR$ is locally Lipschitzian around $\ox$; see, e.g., \cite[Theorem~3.86]{m-book1}.
Now we apply inclusion \eqref{eq:levelset} to each set $\O_i$ in \eqref{ine} and substitute this into the NQC
\eqref{eq:QC-N(O)} as well as into the qualification condition \eqref{sin1} and the optimality conditions \eqref{u-fr1}
and \eqref{lower1} for problem \eqref{eq:SIP1} with the constraint sets \eqref{ine}. It follows in this way that the
SQC \eqref{eq:QCgi} and all the relationships \eqref{up-ine}--\eqref{lo-sub} imply the aforementioned conditions of
Theorems~\ref{Thm:NecPg1} and \eqref{Thm:Nec2} in the setting \eqref{sip-ine} under consideration. It shows furthermore
that the SCC \eqref{scc} yields the NCC \eqref{eq:CC} for the sets $\O_i$ in \eqref{ine}, which thus completes the proof
of the theorem. $\h$\vspace*{0.05in}

Now we consider in more detail the case of {\em convex} constraint
functions $\ph_i$ in \eqref{sip-ine}. Note that the validity of
the SQC \eqref{eq:QCgi} is ensured in the case by the
interior-type condition \eqref{eq:QCcvex} of
Proposition~\ref{Prop:QCcvex}. The next theorem justifies
necessary optimality conditions for problems with countable convex
inequalities, which does not require either interiority-type or
SQC constraint qualifications while containing a qualification
condition that implies both the CHIP and SCC in \eqref{scc}. Let
us first recall this condition; see \cite{DMN09b,DMN09a} and the
references therein. We sat that the SIP problem \eqref{sip-ine}
with the constraints given by convex functions $\ph_i$, $i\in\N$,
satisfies the {\em Farkas-Minkowski constraint qualification}
(FMCQ) if the set
\begin{equation}\label{eq:FM}
\co\Big[\cone\bigcup_{i=1}^\infty\epi\ph_i^*\Big]\;\mbox{ is closed in }\;\R^n\times\R,
\end{equation}
where $\vartheta^*(x^*):=\sup\{\la
x^*,x\ra-\vartheta(x)|\;x\in\R^n\}$ stands for the conjugate
function to $\vartheta\colon\R^n\to\oR$.

\begin{Theorem}\label{sip-convex} {\bf (upper and lower subdifferential conditions for SIP with convex inequality
constraints).} Let all the general assumptions but SQC
\eqref{eq:QCgi} of Theorem~{\rm\ref{Thm:NecPf}} be fulfilled at
the local optimal solution $\ox$ to \eqref{sip-ine}. Assume also
that the constraint functions $\ph_i$, $i\in\N$, are convex. The
both assertion {\rm (i)} and {\rm (ii)} of
Theorem~{\rm\ref{Thm:NecPf}} are satisfied. Furthermore, the
fulfillment of the FMCQ \eqref{eq:FM} implies that the CHIP
\eqref{chip} holds automatically and that the closure operation in
\eqref{up-ine}--\eqref{lo-sub} can be omitted.
\end{Theorem}
{\bf Proof.} Note first that inclusion \eqref{eq:levelset} holds as equality for convex functions, i.e.,
\begin{equation}\label{convex-lev}
N(\ox;\O_i)=\R_+\partial\ph_i(\ox)\;\mbox{ for
}\;\O_i=\{x\in\R^n\big|\;\ph_i(x)\le 0\big\},\quad i\in\N.
\end{equation}
Combining \eqref{convex-lev} with Theorem~\ref{Thm:CHIPchar} and
taking into account that $N(\ox;\O_i)=\{0\}$ when $\ph_i(\ox)<0$,
we can equivalently rewrite the assumed CHIP in the form
\begin{equation}\label{chip-conv}
N\Big(\ox;\bigcap_{i=1}^\infty\O_i\Big)=\cl\co\bigcup_{i\in
J(\ox)}\big[\R_+\partial\ph_i(\ox)\big]\; \mbox{ with
}\;J(\ox):=\big\{i\in\N\big|\;\ph_i(\ox)=0\big\}.
\end{equation}
Substituting the latter into the upper and lower subdifferential optimality conditions
\begin{equation*}
-\Hat\partial^+\ph(\ox)\subset
N\Big(\ox;\bigcap_{i=1}^\infty\O_i\Big)\;\mbox{ and
}\;0\in\partial \ph(\ox)+N\Big(\ox;\bigcap_{i=1}^\infty\O_i\Big)
\end{equation*}
for problem \eqref{sip-ine}, which follow from
\cite[Prositions~5.2 and 5.3]{m-book2}, respectively, we arrive at
the conclusions in (i) and (ii) of Theorem~\ref{Thm:NecPf}.

To complete the proof of the theorem, it remains to check that the
FMCQ \eqref{eq:FM} simultaneously implies the fulfillments of the
CHIP \eqref{chip} and the SCC \eqref{scc}. It follows from
\cite[Corollary~3.6]{DMN09b} that the FMCQ yields the
representation
\begin{equation}\label{dmn}
N\Big(\ox,\bigcap_{i=1}^\infty\O_i\Big)=\bigcup_{\lm\in A(\ox)}\Big[\sum_{i\in J(\ox)}\lm_i\partial\ph_i(\ox)\Big]
\end{equation}
for the constraint sets $\O_i$, where $A(\ox)$ denotes the
collection of Lagrange multipliers $\lm=(\lm_i)_{i\in\N}$ such
that $\lm\in A(\ox)$ if and only if $\lm_i\ge 0$ for $i\in J(\ox)$
and $\lm_i=0$ otherwise. We obviously have from \eqref{convex-lev}
and \eqref{dmn} that
\begin{equation}\label{conv1}
N\Big(\ox,\bigcap_{i=1}^\infty\O_i\Big)=\co\bigcup_{i\in J(\ox)}\R_+\partial\ph_i(\ox)=\co\bigcup_{i=1}^\infty
N(\ox;\O_i).
\end{equation}
Since the normal cone $N(\ox;\O)$ is closed, it follows from
\eqref{conv1} that the set co$\big\{\bigcup_{i\in
J(\ox)}[\R_+\partial\ph_i(\ox)]\big\}$ is closed as well; the
latter is clearly equivalent to the SCC \eqref{scc} at $\ox$. On
the other hand, we have from \eqref{conv1} that the strong CHIP
\eqref{eq:sCHIP1} holds, which implies the fulfillment of the CHIP
\eqref{chip} by Theorem~\ref{Thm:CHIPchar} and thus completes the
proof of this theorem. $\h$\vspace*{0.05in}

Next we present efficient specifications of both upper and lower
subdifferential optimality conditions from
Theorem~\ref{sip-convex} for SIP with linear inequality
constraints. In the finite-dimensional countable case under
consideration the results obtained in this way reduce to those
from \cite[Theorems~3.1 and 4.1]{CLMParra09b} while it is {\em
not} assumed here the strong Slater condition and the coefficient
boundedness imposed in \cite{CLMParra09b}. For simplicity we
consider the case of homogeneous constraints and suppose that
$\ox=0$ is a local optimal solution.

\begin{Corollary}\label{sip-lin} {\bf (upper and lower subdifferential conditions for SIP with linear inequality
constraints).} Let $\ox=0$ be a local optimal optimal solution to the SIP problem
\begin{equation}\label{sip-lin1}
\mbox{minimize }\;\ph(x)\;\mbox{ subject to }\;\la a_i,x\ra\le 0\;\mbox{ for all }\;i\in\N,
\end{equation}
where $\ph\colon\R^n\to\oR$ is finite at the origin. Then we have the inclusions
\begin{equation}\label{up-lin}
-\Hat\partial^+\ph(0)\subset\cl\co\Big[\bigcup_{i=1}^\infty\big\{\lm a_i\big|\;\lm\ge 0\big\}\Big].
\end{equation}
\begin{equation}\label{lo-lin}
0\in\partial\ph(0)+\cl\co\Big[\bigcup_{i=1}^\infty\big\{\lm a_i\big|\;\lm\ge 0\big\}\Big],
\end{equation}
where \eqref{lo-lin} holds provided that $\ph$ is lower semicontinuous around the origin and
\begin{equation}\label{qc-lin}
\Big(\cl\co\Big[\bigcup_{i=1}^\infty\big\{\lm a_i\big|\;\lm\ge 0\big\}\Big]\Big)\cap\big(-\partial^
\infty\ph(0)\big)=\{0\}.
\end{equation}
Furthermore, the FMCQ implies that the closure operations can be omitted in \eqref{up-lin}--\eqref{qc-lin}.
\end{Corollary}
{\bf Proof}. Since the CHIP is automatic for the linear inequality
system in \eqref{sip-lin1} at the origin and by
Corollary~\ref{linear-SIP} we have the normal cone representation
\eqref{lin1}, all the results of this corollary follow from the
corresponding results of Theorem~\ref{sip-convex}.
$\h$\vspace*{0.05in}

Finally in this section, we present several examples illustrating
the qualification conditions imposed in Theorem~\ref{sip-convex}
and their comparison with known results in the in the literature.

\begin{Example}\label{ex-convex} {\bf (comparison of qualification conditions).}
{\rm All the examples below concern lower subdifferential
conditions for SIP problems \eqref{sip-ine} with convex cost and
constraint functions.

{\bf (i)} {\bf The CHIP \eqref{chip} and the SCC \eqref{scc} are
independent}. Consider a linear constraint system in
\eqref{sip-lin} at $\ox=0)\in\R^2$ for $\ph_i(x)=\la a_i,x\ra$
with $a_i=(1,i)$ as $i\in\N$, which has the CHIP. At the same time
the set
\begin{equation*}
\co\bigcup_{i=0}^\infty\R_+\partial\ph_i(\ox)=
\co\big\{\lm(1,i)\in\R^2\big|\;\lm\ge 0,\;i\in\N\big\}=\R_+^2\setminus\big\{(0,\lm)\big|\;\lm>0\big\}
\end{equation*}
is not closed, and hence the SCC \eqref{scc} does not hold. On the
other hand, for the quadratic functions $\ph_i(x)=i x^2_i-x_2$ as
$i\in\N$ as $x=(x_1,x_2)\in\R^2$, we get
$\partial\ph_i(0)=\nabla\ph_i(0)=(0,-1)$, and hence the SCC
\eqref{scc} holds at the origin while the CHIP is violated at this
point by Example~\ref{Ex:CHIPfail1}(ii).

{\bf(ii)} {\bf (CHIP and SCC versus FMCQ and CQC).} Besides the
FMCQ \eqref{eq:FM}, another qualification condition is employed in
\cite{DMN09b,DMN09a} to obtain necessary optimality conditions of
{\em Karush-Kuhn-Tucker} (KKT) type (no closure operation
in\eqref{lo-sub}) for {\em fully convex} SIP problems
\eqref{sip-ine} involving all the convex functions $\ph$ and
$\ph_i$. This condition, named the {\em closedness qualification
condition} (CQC) is formulated as follows via the convex conjugate
functions: the set
\begin{equation}\label{eq:CQC}
\epi\ph^*+\co\Big[\cone\bigcup_{i=1}^\infty\epi\ph_i^*\Big]\;
\mbox{ is closed in }\;\R^n\times\R.
\end{equation}
It is obvious that the FMCQ implies the CQC while the latter is
implied only for fully convex SIP problems. The next example
presents a fully convex SIP problem satisfying both CHIP and SCC
but not the CQC (and hence not FMCQ). This shows that
Theorem~\ref{sip-convex} holds in this case to produce the KKT
optimality condition while the corresponding result of
\cite{DMN09b} is not applicable.

Consider the SIP \eqref{sip-convex} with $x=(x_1,x_2)\in\R^2$, $\ox=(0,0)$, $\ph(x)=-x_2$, and
\begin{equation*}
\ph_i(x_1,x_2)=\left\{\begin{array}{ll}
ix_1^2-x_2&\mbox{ if }\;x_1<0,\\
-x_2&\mbox{ if }\; x_1\ge 0,
\end{array}\right.
\quad i\in\N.
\end{equation*}
We have $\partial\ph_i(\ox)=\nabla\ph_i(\ox)=(0,-1)$ for all $i\in\N$,
and hence the SCC (\ref{scc}) holds. It is easy to check that the CHIP holds at $\ox$, since
\begin{equation*}
T\Big(\ox;\bigcap_{i=1}^\infty\O_i\Big)=T(\ox;\O_i)=\R\times\R_+\;\mbox{
for }\;\O_i:=\big\{x\in\R^2\big|\; \ph_i(x)\le 0\big\},\quad
i\in\N.
\end{equation*}
On the other hand, for $x^*=(\lm_1,\lm_2)\in\R^n$ we compute the conjugate functions by
\begin{equation*}
\ph^*(x^*)=\left\{\begin{array}{ll}
0&\mbox{ if }\;(\lm_1,\lm_2)=(0,-1),\\
\infty&\mbox{ otherwise }
\end{array}\right.\quad
\mbox{and }\;\ph^*_i(x^*)=\left\{
\begin{array}{ll}
\disp\frac{\lm_1^2}{4i}&\mbox{ if }\;\lm_1\le 0,\;\lm_2=-1,\\
\infty&\mbox{ otherwise}.
\end{array}\right.
\end{equation*}
This shows that the convex sets
\begin{equation*}
\co\Big[\cone\bigcup_{i=0}^\infty\epi\ph_i^*\Big]\quad\mbox{and}\quad\epi
\ph^*+\co\Big[\cone\bigcup_{i=0}^\infty\epi\ph_i^*\Big]
\end{equation*}
are not closed in $\R^2\times\R$, and hence the FMCQ \eqref{eq:FM} and the CQC \eqref{eq:CQC} are not satisfied.}
\end{Example}

\section{Applications to Multiobjective Optimization}
\setcounter{equation}{0}

The last section of this paper concerns problems of multiobjective
optimization with set-valued objectives and countable constraints.
Although optimization problems with single-valued/vector and (to a
lesser extent) set-valued objectives have been widely considered
in optimization and equilibrium theories as well as in their
numerous applications (see, e.g., the books
\cite{grtz,Jahn-book04,m-book2} and the references therein), we
are not familiar with the study of such problems involving
countable constraints. Our interest is devoted to deriving
necessary optimality conditions for problems of this type based on
the dual-space approach to the general multiobjective optimization
theory developed in \cite{B-M10a,B-M10b,m-book2} and the new
tangential extremal principle established in \cite{MorPh10a}.

The main problem of our consideration is as follows:
\begin{equation}\label{eq:Mob}
\mbox{minimize }\;F(x)\;\mbox{ subject to }\;
x\in\O:=\bigcap_{i=1}^\infty\O_i\subset\R^n,
\end{equation}
where $\O_i$, $i\in\N$, are closed subsets of $\R^n$, where
$F\colon\R^n\tto\R^m$ is a set-valued mapping of closed graph, and
where ``minimization" is understood with respect to some partial
ordering ``$\le$'' on $\R^m$. We pay the main attention to the
multiobjective problems with the {\em Pareto-type} ordering:
\begin{equation*}
y_1\le y_2\ \mbox{ if and only if }\ y_2-y_1\in\Theta,
\end{equation*}
where $\emp\ne\Theta\subset\R^m\setminus\{0\}$ is a closed,
convex, and pointed ordering cone. In the aforementioned
references the reader can find more discussions on this and other
ordering relations.

Recall that a point $(\ox,\oy)\in\gph F$ with $\ox\in\O$ is a {\em
local minimizer} of problem (\ref{eq:Mob}) if there exists a
neighborhood $U$ of $\ox$ such that there is no $y\in F(\O\cap U)$
preferred to $\oy$, i.e.,
\begin{equation}\label{eq:minMob}
F(\O\cap U)\cap(\oy-\Theta)=\{\oy\}.
\end{equation}
Note that notion \eqref{eq:minMob} does not take into account the
image localization of minimizers around $\oy\in F(\ox)$, which is
essential for certain applications of set-valued minimization,
e.g., to economic modeling; see \cite{B-M10b}. A more appropriate
notion for such problems is defined in \cite{B-M10b} under the
name of {\em fully localized minimizers} as follows: there are
neighborhoods $U$ of $\ox$ and $V$ of $\oy$ such that
\begin{equation}\label{eq:LminMob}
F(\O\cap U)\cap(\oy-\Theta)\cap V=\{\oy\}.
\end{equation}

The next result establishes necessary optimality conditions of the
coderivative type for fully localized minimizers of problem
\eqref{eq:Mob} with countable constraints based on the approach of
\cite{m-book2} to problems of multiobjective optimizations, its
implementations in \cite{B-M10a,B-M10b} specifically for problems
with set-valued criteria, and the tangential extremal principle
for countable sets \cite{MorPh10a}. We address here fully
localized minimizers for multiobjective problems \eqref{eq:Mob}
with normally regular feasible sets, i.e., when $N(\ox;\O)=\Hat
N(\ox;\O)$, which particularly includes the case of convex set
$\O_i$, $i\in\N$.

\begin{Theorem}\label{nc-cod} {\bf (optimality conditions for
fully localized minimizers of multiobjective problems with
countable constraints and normally regular feasible sets).} Let
the pair $(\ox,\oy)\in\gph F$ be a fully localized minimizer for
\eqref{eq:Mob} with the CHIP system of countable constraints
$\{\O_i\}_{i\in\N}$. Assume that the feasible set
$\O=\bigcap_{i=1}^\infty\O_i$ is normally regular at $\ox\in\O$
and that the NQC \eqref{eq:QC-N(O)} and the coderivative
qualification condition
\begin{equation}\label{qc-mo}
D^*F(\ox,\oy)(0)\bigcap\Big[-\cl\Big\{\sum_{i\in
I}x^*_i\Big|\;x^*_i\in N(\ox;\O_i),\;I\in{\cal
L}\Big\}\Big]=\big\{0\big\}
\end{equation}
are satisfied. Then there is $0\ne y^*\in-N(0;\Theta)$ such that
\begin{equation}\label{nc-mo}
0\in D^*F(\ox,\oz)(y^*)+\cl\Big\{\sum_{i\in I}x^*_i\Big|\;
x^*_i\in N(\ox;\O_i),\;I\in{\cal L}\Big\}.
\end{equation}
\end{Theorem}
{\bf Proof.} Applying \cite[Theorem~3.4]{B-M10b} for fully
localized minimizers of set-valued optimization problems with
abstract geometric constraints $x\in\O$ (cf.\ also
\cite[Theorem~5.3]{B-M10a} for the case of local minimizers
\eqref{eq:minMob} and \cite[Theorem~5.59]{m-book2} for vector
single-objective counterparts), we find
\begin{equation}\label{nc-mo1}
0\ne-y^*\in N(0;\Th)\;\mbox{ and }\;x^*\in
D^*F(\ox,\oy)(y^*)\cap\big(-N(\ox;\O)\big)
\end{equation}
provided the fulfillment of the qualification condition
\begin{equation}\label{qc-mo1}
D^*F(\ox,\oy)(0)\cap\big(-N(\ox;\O)\big)=\{0\}.
\end{equation}
To complete the proof of the theorem, it suffices to employ in
\eqref{nc-mo1} and \eqref{qc-mo1} the sum rule for countable set
intersections from Theorem~\ref{Thm:FnorO} by taking into account
the assumed normal regularity of the intersection set $\O$ at
$\ox$. $\h$\vspace*{0.05in}

Note that the qualification condition \eqref{qc-mo} holds
automatically if the objective mapping $F$ is {\em Lipschitz-like}
(or has the Aubin property) around $(\ox,\oy)\in\gph F$, i.e.,
there are neighborhoods $U$ of $\ox$ and $V$ of $\oy$ such that
\begin{equation*}
F(x)\cap V\subset F(u)+\ell\|x-u\|\B\;\mbox{ for all }\;x,u\in U
\end{equation*}
with some number $\ell\ge 0$. Indeed, it follows from the
Mordukhovich criterion in \cite[Theorem~9.40]{Rockafellar-Wets-VA}
(see also \cite[Theorem~4.10]{m-book1} and the references therein)
that $D^*F(\ox,\oy)(0)=\{0\}$ in this case.\vspace*{0.05in}

Next we introduce two kinds of ``graphical" minimizers for
multiobjective problems for which, in particular, we can avoid the
normal regularity assumption in optimality conditions of type
\eqref{nc-mo} in Theorem~\ref{nc-cod}. The definition below
concerns multiobjective optimization problems with general
geometric constraints that may not be represented as countable set
intersections.

\begin{Definition}\label{Def:GTmin}{\bf (graphical and tangential
graphical minimizers).} Let $(\ox,\oy)\in\gph F$ with $\ox\in\O$.
We say that:

{\bf (i)} $(\ox,\oy)$ is a {\sc local graphical minimizer} to
problem \eqref{eq:Mob} if there are neighborhoods $U$ of $\ox$ and
$V$ of $\oy$ such that
\begin{equation}\label{eq:Gmin}
\gph F\cap\Big[\O\times(\oy-\Theta)\Big]\cap\big(U\times
V\big)=\big\{(\ox,\oy)\big\}.
\end{equation}

{\bf (ii)} $(\ox,\oy)$ is a {\sc local tangential graphical
minimizer} to problem \eqref{eq:Mob} if
\begin{equation}\label{eq:Tmin}
T\big((\ox,\oy);\gph F\big)\cap
\big[T(\ox;\O)\times(-\Theta)\big]=\{0\}.
\end{equation}
\end{Definition}

Similarly to the discussions and examples on relationships between
local extremal and tangentially extremal points of set systems
given in \cite{MorPh10a}, we observe that the optimality notions
in Definition~\ref{Def:GTmin} are {\em independent} of each other.
Let us now compare the the graphical optimality of
Definition~\ref{Def:GTmin}(i) with fully localized minimizers of
\eqref{eq:LminMob}.

\begin{Proposition}\label{graph} {\bf (relationships between fully localized and
graphical minimizers).} Let $(\ox,\oy)\in\gph F$ be a feasible
solution to problem \eqref{eq:Mob} with general geometric
constraints. Then the following assertions are satisfied:

{\bf (i)} $(\ox,\oy)$ is a local graphical minimizer if it is a
fully localized minimizer for this problem.

{\bf (ii)} The opposite implication holds if there is a
neighborhood $U$ of $\ox$ such that $\oy\notin F(x)$ for every
$\ox\ne x\in\O\cap U$.
\end{Proposition}
{\bf Proof}. To justify (i), assume that $(\ox,\oy)$ is a local
graphical minimizer, take its neighborhood $U\times V$ from
Definition~\ref{Def:GTmin}(i), and pick any
\begin{equation*}
y\in F(\O\cap U)\cap(\oy-\Theta)\cap V.
\end{equation*}
Then there is $x\in\O\cap U$ such that $y\in F(x)$, and so
\begin{equation*}
(x,y)\in \gph F\cap\big[\O\times(\oy-\Theta)\big]\cap\big(U\times
V\big)=\big\{(\ox,\oz)\big\}.
\end{equation*}
Thus $F(\O\cap U)\cap(\oy-\Theta)\cap V=\{\oy\}$, i.e.,
$(\ox,\oy)$ is a fully localized minimizer for \eqref{eq:Mob}.

Next we prove (ii). Suppose that $(\ox,\oy)$ is a fully localized
minimizer with a neighborhood $U\times V$, shrink $U$ so that the
assumption in (ii) holds, and take
\begin{equation*}
(x,y)\in\gph F\cap\big[\O\times(\oy-\Theta)\big]\cap\big(U\times
V\big).
\end{equation*}
Since $y\in F(x)$, it follows that $y\in F(\O\cap U)\cap
(\oy-\Theta)\cap V=\{\oy\}$. If $x\ne\ox$, the latter contradicts
the assumption in (ii). Thus $x=\ox$, which completes the proof of
the proposition. $\h$\vspace*{0.05in}

The next theorem uses the full strength of the tangential extremal
principle of \cite{MorPh10a} justifying the necessary optimality
conditions of Theorem~\ref{nc-cod} for {\em tangential} graphical
minimizers of the multiobjective problem \eqref{eq:Mob} with
countable constraints without imposing the normal regularity
requirement of the feasible set.

\begin{Theorem}\label{Thm:Nec-Tmin} {\bf (optimality conditions
for tangential graphical minimizers).} Let $(\ox,\oy)$ be a local
tangential graphical minimizer for problem \eqref{eq:Mob} under
the fulfillment all the assumptions of Theorem~{\rm\ref{nc-cod}}
but the normal regularity of $\O$ at $\ox$. Suppose in addition
that $\rm{int}\Th\ne\emp$. Then there is $0\ne y^*\in-N(0;\Theta)$
such that the necessary optimality condition \eqref{nc-mo} is
satisfied.
\end{Theorem}
{\bf Proof}. We have by Definition~\ref{Def:GTmin}(ii) that $\disp
T((\ox,\oz);\gph F)\cap\big[\Lm\times(-\Theta)\big]=\{0\}$ with
$\Lm:=T(\ox;\O)$. Since the system $\{\O_i\}_{i\in\N}$ has the
CHIP at $\ox$, it follows that
\begin{equation*}
\Lm=\bigcap_{i=1}^\infty\Lm_i\;\mbox{ with }\;\Lm_i:=T(\ox;\O_i).
\end{equation*}
Further, define the closed cones $\GG_0:=T((\ox,\oz);\gph F)$ and
$\GG_i:=\Lm_i\times(-\Theta)$ as $i\in\N$ with
$\displaystyle\bigcap_{i=0}^\infty \GG_i=\{0\}$ and show that for
any $\xi\in\Theta$ we get
\begin{equation}\label{tan3}
\bigcap_{i=1}^\infty\GG_i\bigcap\Big[\GG_0+(0,\xi)\Big]=\emp.
\end{equation}
Indeed, supposing the contrary gives us a vector
$(x,y)\in\R^n\times\R^m$ with $(x,y-\xi)\in\GG_0$ and
$(x,y)\in\Lm_i\times(-\Theta)$ for all $i\in\N$. Since $\Theta$ is
a closed and convex cone, we also have the inclusion $(x,y-\xi)\in
\Lm_i\times(-\Theta)=\GG_i$ as $i\in\N$, and hence
\begin{equation*}
(x,y-\xi)\in\bigcap_{i=0}^\infty\GG_i=\{0\}.
\end{equation*}
It follows therefore that $y=\xi\in-\Theta$, which implies by the
pointedness of the cone $\Theta$ that
$\xi\in(-\Theta)\cap\Theta=\{0\}$, a contradiction justifying
\eqref{tan3}.

The latter means that $\{\GG_i\},i=0,1,\ldots$, is a countable
system of cones extremal at the origin with the nonoverlapping
condition $\bigcap_{i=0}^\infty\GG_i=\{0\}$. Now applying the
tangential extremal principle of Theorem~\ref{Thm:TEPinf} to this
system of cones and using also \cite[Proposition~2.1]{MorPh10a},
we get elements $(x^*_i,y^*_i)$ as $i=0,1,\ldots$ satisfying the
relationships
\begin{equation}\label{tan4}
(x^*_0,y^*_0)\in N(0;\GG_0)\subset N\big((\ox,\oy);\gph F\big),
\end{equation}
\begin{equation}\label{tan5}
(x^*_i,y^*_i)\in N(0;\GG_i)\subset
N(\ox;\O_i)\times\big[-N(0;\Theta)\big],\quad i\in\N,
\end{equation}
\begin{equation}\label{tan6}
\sum_{i=0}^\infty
\frac{1}{2^i}\Big(x^*_i,y^*_i\Big)=0,\quad\mbox{and}\quad
\sum_{i=0}^\infty\frac{1}{2^i}\Big(\|x^*_i\|^2+\|y^*_i\|^2\Big)=1.
\end{equation}
It follows from \eqref{tan4}--\eqref{tan6} that
\begin{equation}\label{tan7}
x^*_0\in D^*F(\ox,\oy)(-y^*_0)\;\mbox{ and }\;
\disp-y^*_0=\sum_{i=1}^\infty\frac{1}{2^i}y^*_i\in-N(0;\Theta),
\end{equation}
where the latter inclusion holds by the convexity and closedness
of the cone $N(0;\Theta)$.

There are the two possible cases in \eqref{tan7}: $y^*_0\ne 0$ and
$y^*_0=0$. In the first case we get
\begin{equation*}
0\in D^*F(\ox,\oy)(-y^*_0)+\sum_{i=1}^\infty\frac{1}{2^i}x^*_i,
\end{equation*}
which readily implies the optimality condition \eqref{nc-mo} with
$0\ne y^*:=-y^*_0\in-N(0;\Th)$; cf.\ the proof of the second part
of \cite[Theorem~5.4]{MorPh10a}.

To complete the proof of this theorem, it remains to show that the
case of $y^*_0=0$ in \eqref{tan7} cannot be realized under the
imposed qualification conditions \eqref{eq:QC-N(O)} and
\eqref{qc-mo}. Indeed, for $y^*_0=0$ we have from \eqref{tan5} and
\eqref{tan7} that
\begin{equation}\label{tan8}
-\frac{1}{2}y^*_1=\sum_{i=2}^\infty\frac{1}{2^i}y^*_i\in
\big[-N(0;\Theta)\big]\cap N(0;\Theta).
\end{equation}
Since the cone $\Th$ is convex, it follows from \eqref{tan8} that
\begin{equation*}
\la y^*_1,y\ra\le 0\;\mbox{ and }\;\la y^*_1,y\ra\ge 0\;\mbox{ for
any }\;y\in\Th,
\end{equation*}
i.e., $\la y^*_1,y\ra=0$ on $\Th$. The latter implies that
$y^*_1=0$ by $\rm{int}\Th\ne\emp$.

Proceeding in this way by induction gives us that $y^*_i=0$ for
all $i\in\N$. Now it follows from \eqref{tan5} and the first
inclusion in \eqref{tan7} that $x^*_0=0$ by the assumed
coderivative qualification condition \eqref{qc-mo}. Hence we get
from \eqref{tan6} the relationships
\begin{equation*}
\sum_{i=0}^\infty\frac{1}{2^i}x^*_i=0\;\mbox{ and
}\;\sum_{i=0}^\infty\frac{1}{2^i}\|x^*_i\|^2=1,
\end{equation*}
which contradict the assumed NQC \eqref{eq:QC-N(O)} and thus
complete the proof of the theorem. $\h$\vspace*{0.05in}

Note in conclusion that, similarly to Section~4, we can develop
necessary optimality conditions for multiobjective problems with
countable constraints of operation and inequality types.

\end{document}